\documentclass[12pt]{report}
\usepackage[centertags]{amsmath}
\usepackage{amsfonts}
\usepackage{amssymb}
\usepackage{amsthm}
\usepackage{newlfont}
\usepackage{verbatim}
\usepackage[dvips]{graphics}
\usepackage{enumerate}
\newtheorem{proposition}{Proposition}[section]
\newtheorem{lemma}[proposition]{Lemma}
\newtheorem{definition}[proposition]{Definition}

\newtheorem{corollary}[proposition]{Corollary}
\newtheorem{theorem}[proposition]{Theorem}

\newtheorem{definition_ch}{Definition}[chapter]
\newtheorem{corollary_ch}[definition_ch]{Corollary}
\newtheorem{theorem_ch}[definition_ch]{Theorem}

\newcommand{\me}{\mathrm{e}}
\newcommand{\mi}{\mathrm{i}}
\newcommand{\dif}{\mathrm{d}}

\newcommand{\D}{\ensuremath{\mathbb{D}}}
\newcommand{\Dstar}{\ensuremath{\mathbb{D}^*}}
\newcommand{\SOG}{\ensuremath{S^O(\Gamma)}}
\newcommand{\SCG}{\ensuremath{S^C(\Gamma)}}
\newcommand{\dd}[1]{{\textstyle\frac{\partial}{\partial #1}}}

\begin{document}
  \begin{titlepage}
\begin{center}
 {\bfseries\textbf \large RIEMANN SURFACES AND
 3-REGULAR GRAPHS}
\vfill
Research Thesis
\\[5Ex]
submitted in partial fulfillment of the\\
requirements for the degree of\\
Master of Science in Mathematics
\\[10Ex]
Dan Mangoubi
\vfill
Submitted to the Senate of\\
The Technion - Israel Institute of Technology
\\[5Ex]

Tevet, 5761 \hspace{3Ex} Haifa \hspace{3Ex} January 2001
\end{center}
\end{titlepage}

%----------------
  \begin{titlepage}

\noindent This research thesis was done under the supervision
of Professor Robert Brooks in the Department of Mathematics.
\vfill
\noindent I would like to thank Professor Brooks for lighting my way,
caring and loving.

\vfill
\noindent The thesis is dedicated to my parents, sister and brother.
Cette th\`{e}se est \'{e}galement dedi\'{e}e \`{a} mes grand-m\`{e}res
bien aim\'{e}es Ir\`{e}ne et Sarah.
\\[18Ex]

\vfill
\noindent The generous financial help of 
the Forchheimer Foundation Fellowship is gratefully acknowledged.

\end{titlepage}

%-----------------
  \addtocontents{toc}{\protect\thispagestyle{empty}}
  \tableofcontents
%-----------------
  \listoffigures
  \thispagestyle{empty}
  \setcounter{page}{0}
%------------------
  \chapter*{Abstract}
\addcontentsline{toc}{chapter}{Abstract}
In this thesis, we consider a way to
construct a rich family of compact Riemann surfaces
in a combinatorial way. Given a 3-regular graph with orientation, 
we construct a finite-area hyperbolic Riemann surface by gluing 
triangles together according to the combinatorics of the graph.
We then compactify this surface by adding finitely many points.

We discuss this construction by considering a number of examples.
In particular, we see that the surface depends in a strong way on the
orientation.

We then consider the effect the process of compactification has on
the hyperbolic metric of the surface. To that end, we ask 
when we can change the metric in the horocycle neighbourhoods of the cusps
to get a hyperbolic metric on the compactification.
In general, the process of compactification can have drastic effects 
on the hyperbolic structure. 
For instance, if we compactify the 3-punctured sphere
we lose its hyperbolic structure.

We show that when the cusps have lengths bigger than $2\pi$ we can 
fill in the horocycle neighbourhoods and retain negative curvature.

Furthermore, the condition that the horocycles
have length $>$ $2\pi$ is sharp. We show by examples that there exist
curves arbitrarily close to horocycles of length $2\pi$, 
which cannot be so filled in. 
Such curves can even be taken to be convex.

  \chapter*{List of Symbols}
\addcontentsline{toc}{chapter}{List of Symbols}
\begin{tabular}{ll}
$\Gamma$             & A (3-regular) graph.\\
$\SOG$               & A non-compact Riemann Surface constructed from 
                        $\Gamma$.\\
$\SCG$               & A compact Riemann Surface constructed from $\Gamma$.\\
$G$                  & A group.\\
$ISO^+$              & Orientation-Preserving Isometries.\\
$\D$                 & $\left\{z\in\mathbb{C}: |z|<1\right\}$.\\
$\Dstar$             & $\D\setminus\{0\}$. \\
$\mathcal{B}_r$      & $\left\{z\in\Dstar:|z|<r\right\}.$\\
$\dif s^2_{\D}$      & The complete hyperbolic metric on $\D$.\\
$\dif s^2_{\Dstar}$  & The complete hyperbolic metric on $\Dstar$.\\
$\kappa(\dif s^2)$   & The curvature of the Riemannian metric $\dif s^2$.\\
$\kappa_g$           & Geodesic curvature.\\
$\Delta$             & The Laplacian.
\end{tabular}

  \chapter*{Introduction}
\addcontentsline{toc}{chapter}{Introduction}
%------------------------------
\section*{Outline of the thesis}
\addcontentsline{toc}{section}{Outline of the thesis}
%------------------------------
This thesis considers a way to construct Riemann surfaces
out of 3-regular graphs.
In Chapter~\ref{ch:construction},
we present a method to construct a finite-area Riemann surface 
out of a 3-regular graph.
We then compactify it to get a compact Riemann surface.
In sections~\ref{sec:ideal_triangles}--~\ref{sec:genus}
we describe this construction in detail.
In section~\ref{sec:automorphisms} we show
that symmetries of the graph are reflected in symmetries
of the surface. This allows us to give several examples in
section~\ref{sec:examples}.
In the last section of the chapter we show that the compact
Riemann surfaces we construct are dense in the moduli space 
by considering them as Belyi surfaces.

We would like to learn about the compact Riemann surfaces
by studying the geometry and spectral properties of the
graphs from which they originated.
Early works of Prof.\ Brooks (\cite{brooks:tc, brooks:vdrm}) 
show connections between the 
graphs and the non-compact Riemann surface.
Later results in~\cite{brooks:platonic} show that under
nice conditions (the ``large cusps condition''), 
the complete metrics of constant curvature
on the non-compact Riemann surface and its compactification
are arbitrarily close, outside of standard cusp neighbourhoods.

It is at this point that we start chapter~\ref{ch:relations}.
We would like to know, for example, when we can control
the metric on the compact Riemann surface. More
precisely, when can we change the hyperbolic metric 
on the non-compact Riemann surface in neighbourhoods of the cusps
to get a hyperbolic metric on the compact Riemann surface?
For example, if we compactify the 3-punctured sphere
we lose the hyperbolic structure. 
To that end, we work in a punctured-disk neighbourhood of a cusp,
and consider the problem of extension of metrics smoothly across
the cusps to obtain metrics of negative curvature.
We prove that in a ``large cusps'' condition we
can in fact extend the metric from the outside to a smooth metric of
negative curvature inside.
This is the content of theorem~\ref{thm:weak_control}.

The ``large cusps'' condition is discussed in section~\ref{sec:large_cusps}.
There, we prove also the classical Shimizu--Leutbecher theorem in a nice way. 
We study the hyperbolic punctured disk in 
sections~\ref{sec:dstar}--\ref{sec:curvature}, 
and theorem~\ref{thm:weak_control} is proved in sections~\ref{sec:key_lemma}
and~\ref{sec:proof}.

In section~\ref{sec:discussion} we show that the ``large cusps''
condition in theorem~\ref{thm:weak_control} is sharp.
This leads us in section~\ref{sec:extension}
to a very interesting problem:
Suppose we have a closed curve bounding a domain $V$ in the unit disk,
and we have a conformal Riemannian metric of negative curvature outside
of $V$. We ask whether we can extend the metric to $V$ in such a way,
that we are left with a conformal Riemannian metric of negative curvature.
The Gauss--Bonnet theorem imposes a natural restriction on the extension,
but this is not enough. We show that we can find curves which satisfy
the Gauss--Bonnet restriction
for which we cannot extend the metric. We do this by using
the maximum principle for subharmonic functions.
Indeed, we show that such curves can be taken to be convex
and arbitrarily close to the horocycle of length $2\pi$.

We also consider curves which satisfy the
Gauss--Bonnet restriction and for which 
no metric of negative curvature,
conformal or nonconformal, can extend the metric from the outside.

In section~\ref{sec:comparison} we show results from~\cite{brooks:platonic}:
If we have large cusps,
the hyperbolic metric on the compactification
is very close to the original hyperbolic metric, 
outside of cusp neighbourhoods.
This is theorem~\ref{thm:comparison}.
The basic idea here is to use the Ahlfors--Schwarz Lemma, which is proved
and discussed in appendix~\ref{app:ahlfors_schwarz}. 
%--------------------------------------------------------------
\section*{Acknowledgements}
\addcontentsline{toc}{section}{Acknowledgements}
%--------------------------------------------------------------
We would like to thank Curt McMullen for giving us
his permission to use his example in section~\ref{sec:examples}.
We would also like to thank Curt McMullen for
discussions with us concerning the thesis, giving us some new ideas to
think about.

We would like to thank Mikhail Katz for his fruitful ideas concerning
the problem of extension of metrics in section~\ref{sec:extension}.

  \chapter{Riemann Surfaces and 3-regular Graphs}
\label{ch:construction}
%--------------------------------------------
\section{Ideal Triangles}
\label{sec:ideal_triangles}
%--------------------------------------------
Let us denote by $\mathbb{D}$ the unit disk, equipped with
the hyperbolic metric.
\begin{definition}
  An ideal triangle in $\mathbb{D}$ is a triangle
  which has 3 vertices at infinity.
\end{definition}

\noindent \textbf{Remark:} Any two ideal triangles are isometric, 
since we may find a M\"{o}bius transformation, 
which takes one onto the other.

\begin{figure}
   $$\scalebox{0.75}[0.75]{\includegraphics{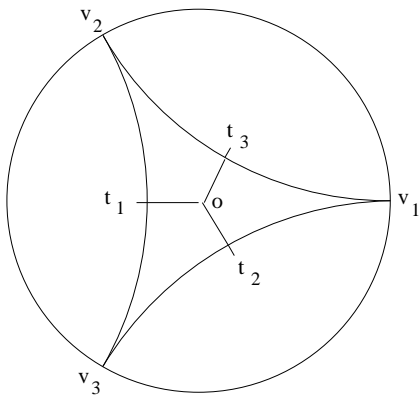}}$$
   \caption{The standard marked ideal triangle}
   \label{fig:triangle}
\end{figure}

Choose a point on each edge of the ideal triangle.
The chosen points will be called \emph{tick-marks}.
\begin{definition}
  A \emph{marked ideal triangle} is an ideal triangle 
  with a tick-mark on each one of its three sides.
\end{definition}
An isomorphism between two marked ideal triangles is
an isomorphism between the ideal triangles which preserves 
the tick-marks.\\
\emph{The standard marked ideal triangle} (figure~\ref{fig:triangle})
is any marked ideal triangle which is isometric 
to the marked ideal triangle whose vertices 
in the disk model are given by
$$ v_1=1,\ v_2=\omega,\ v_3=\omega^2 ,$$ and whose tick-marks are
$$t_1 = -(2-\sqrt{3}),\ 
  t_2 = -(2-\sqrt{3})\omega,\ 
  t_3 = -(2-\sqrt{3})\omega^2,$$
where $\omega=\me^{2\pi \mi/3}$.

\begin{figure}
  $$\scalebox{0.55}[0.55]{\includegraphics{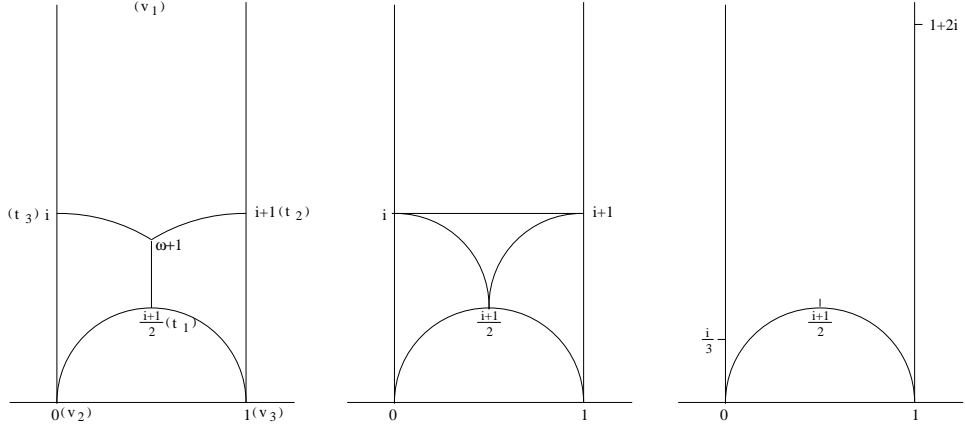}}$$
  \caption{standard and non-standard marked ideal triangles}
  \label{fig:triangle.poincare}
\end{figure}

Under the isometry which takes the half-plane 
  onto the disk and is given by the map
  $$z\mapsto\frac{z-(\omega+1)}{z-(\bar{\omega}+1)}\,$$
  the ideal triangle of figure~\ref{fig:triangle} corresponds to
  the first ideal triangle from the left in 
  figure~\ref{fig:triangle.poincare}.
  In the second ideal triangle there, we have drawn also 
  three horocyclic segments, which will be called
  the \emph{standard horocyclic segments}. 
  Note that each standard horocyclic segment is of length 1.
%-------------------------------------
\section{Some notions in Graph Theory}
\label{sec:graphs}
\begin{figure}
  $$\includegraphics{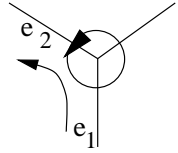}$$
  \caption{A left-hand-turn path}
  \label{fig:lht}
\end{figure}

%-------------------------------------
Throughout our work, the term \emph{graph}
will refer to a graph with loops and multiple edges permitted.
\begin{definition}
A \emph{3-regular graph} is a graph for which there
are three edges incident at each vertex
(loops are counted twice).
\end{definition}
Let $\Gamma$ be a 3-regular graph.
Denote by $V$ and $E$ the set of its vertices and the set of its 
edges respectively.
\begin{definition}
  An \emph{orientation} at a vertex $v$ is a cyclic ordering
  of the edges incident at $v$.
\end{definition}
At each vertex of $\Gamma$ choose an orientation
(We will suppress the orientation and will denote the graph
with orientation by $\Gamma$).
We may regard the orientation as follows:
If we take a tour on the graph, at each vertex
it tells us which way is left.
We define:
\begin{definition}
A \emph{left-hand-turn path} in $\Gamma$
is a directed closed path in $\Gamma$ such that
if $e_1, e_2$ are successive edges in the path meeting at $v$,
then $e_2, e_1$ are successive edges with respect to the
orientation at $v$ (see figure~\ref{fig:lht}).
\end{definition}
Where convenient, we will denote a left-hand-turn path by
a cyclic sequence $\left((v_1,e_1),\ldots,(v_n,e_n)\right)$, where
$(v_i,e_i) \in (V,E)$ and $\partial e_i =\{v_i, v_{i+1}\}$.
In this notation $(v_{i+n},e_{i+n}) = (v_i,e_i)$, and
$(v_i,e_i) \neq (v_j,e_j)$ for $i\not\equiv j \pmod{n}$.
%--------------------------------------------
\section{Triangulations and 3-regular Graphs}
\label{sec:triangles-graphs}
%--------------------------------------------
Let $S$ be a closed (topological) surface.
Let $\mathcal{T}$ be a triangulation of $S$.
Mark a point in the interior of each triangle, and
connect points corresponding to adjacent triangles.
We get a 3-regular graph, $\Gamma$, drawn on the surface.
$\Gamma$ inherits an orientation from the surface's
orientation - at each vertex, $v$, of $\Gamma$ we order the three
emanating edges from $v$ counterclockwise.
For convenience, we will refer to the vertices and edges 
of the triangulation as $\mathcal{T}$-vertices and $\mathcal{T}$-edges 
respectively,
while the vertices and edges of $\Gamma$ will be referred as
$\Gamma$-vertices and $\Gamma$-edges respectively.

First, we note that
any $\Gamma$-edge is crossed by 
exactly one $\mathcal{T}$-edge, and vice versa.
Therefore the $\Gamma$-edges are in one-to-one
correspondence with the $\mathcal{T}$-edges.
Second, any left-hand-turn path in $\Gamma$ 
cuts, on its left side, a disk from the surface
in which there are no $\Gamma$-vertices,
since otherwise one could join any $\Gamma$-vertex
inside the disk by a shortest path in $\Gamma$ 
to the boundary of the disk to get that the boundary
would not be a left-hand-turn path.
In this disk there is exactly one $\mathcal{T}$-vertex:
Indeed, any $\Gamma$-edge which crosses the boundary of the disk
must have one of its ends inside the disk - so
there is at least one $\mathcal{T}$-vertex in this disk.
Conversely, take any $\mathcal{T}$-vertex, $p$, in the disk, 
and let $l$ be the length of the boundary of the disk. 
If there were less than $l$ $\mathcal{T}$-edges emanating from $p$,
one of the triangles surrounding $p$ would contain at least two  
$\Gamma$-vertices from the boundary of the disk.
Therefore, there are $l$ $\mathcal{T}$-edges emanating from
$p$, and we are left with ``no room''
for another $\mathcal{T}$-vertex inside the disk.
%(temp- Did I use orientability of $S$
%only to define an orientation on $\Gamma$? Yes).

We have shown that we have a well-defined bijection:
A left-hand-turn path in $\Gamma$ maps
to the $\mathcal{T}$-vertex it encloses on its left side.
As we will see in the next section,
this bijection lets us reconstruct the surface from the
graph with orientation. Or, in more fancy language,
\{graphs with orientations\} $\cong$
\{triangulations of surfaces\}/\{isotopy\}.
%----------------------------------------------------------------
\section{Constructing a Riemann Surface out of a 3-regular Graph}
\label{sec:construction}
%----------------------------------------------------------------
Let $\Gamma$ be a 3-regular graph with orientation.
Paste one copy of the standard marked ideal triangle on each vertex of
$\Gamma$, in such a way that the three successive edges emanating from
the vertex of the graph fit the three geodesic segments
$ot_1, ot_2, ot_3$ in figure~\ref{fig:triangle}, respectively.
We glue adjacent sides of marked triangles, such that
the tick-marks fit and the orientations
of the corresponding triangles match up.
We can begin doing the pasting without leaving the hyperbolic plane,
until we get a polygon, $\Pi$, together with a side pairing
(Remark: We may think of the graph $\Gamma$ as made out of threads,
and of each triangle as having a white face
and a black face.
The threads are glued on the white face of each triangle.
We build the polygon such that all the triangles have white face up).
Then, we apply the Poincar\'{e} Polygon Theorem (see Appendix~\ref{app:PPT})
to $\Pi$:
We attach to each side pairing $(s_i, s_j)$ an orientation preserving
transformation, $A_{ij}$, such that $A_{ij}(s_i)=s_j$, 
$A_{ij}$ preserves tick-marks,
and $\Pi\cap A_{ij}(\Pi)=\emptyset$.
Denote by $G$ the subgroup of $\mathrm{ISO}^+(\mathbb{D})$
generated by all the side-pairing transformations.
By the Poincar\'{e} Polygon Theorem, $G$ is a discrete group
of isometries with
$\Pi$ as its fundamental domain. Hence $\mathbb{H}^2/G$ is 
a complete hyperbolic Riemann surface, which we will denote by \SOG.

As was shown in the preceding section, 
the cusps of \SOG{} are in bijection with
the left-hand-turn paths in $\Gamma$,
and the triangulation of
\SOG{} obtained by our construction corresponds to $\Gamma$
in the sense of section~\ref{sec:triangles-graphs}.
Finally, take the unique conformal 
compactification (see section~\ref{sec:conformal}) of \SOG,
to get a closed Riemann surface, \SCG.
%---------------------------------------------
\section{A slightly more general construction}
%---------------------------------------------
Again, let $\Gamma$ be a 3-regular graph with orientation.
We paste a marked ideal triangle on each vertex as before,
but now we allow also non-standard marked ideal triangles to be pasted.
The marked ideal triangles can be parametrized by three real
parameters, ($\alpha_1$, $\alpha_2$, $\alpha_3$), which describe
the hyperbolic shift of the tick-marks with respect to the tick-marks
$t_1$, $t_2$, $t_3$ of the standard marked ideal triangle.
For instance, the standard marked ideal triangle will be denoted by
$(0,0,0)$, and the non-standard marked ideal triangle
in figure~\ref{fig:triangle.poincare} is $(0, \ln 2, -\ln 3)$.
In the next step, we glue adjacent edges as before,
remembering to match orientations and tick-marks.

To each pair of a vertex, $v$, and an edge 
$e$ emanating from $v$, there corresponds one side of
the triangle we paste at $v$. Therefore we can
attach to such a pair a real number $\alpha(v,e)$ which describes
the shift of the tick-mark on the corresponding side.
%(temp-loops might be a problem here)

In order to apply the Poincar\'{e} Polygon Theorem
all the vertex-cycle transformations must be parabolic.
This, in turn, is equivalent to the following condition:\\
For every left-hand-turn path in $\Gamma$, 
$l=\left((v_1,e_1),\ldots, (v_n,e_n)\right)$,
$$\sum_{k=1}^n\left[\alpha(v_k,e_k)
   -\alpha(v_k,e_{k-1})\right] = 0\ .$$
%(temp-Here, we have to be more cautious about loops, I have to think
%how to put it nicely).-No.
%--------------------------------------------
\section{The Genus of \SCG}
\label{sec:genus}
%--------------------------------------------
In section~\ref{sec:construction} we obtained from a 3-regular graph with
orientation, $\Gamma$, a closed Riemann surface $\SCG$.
Let us calculate the genus of \SCG{}.
To that end, let $V$, $E$, $F$, be the number of vertices, edges and faces,
respectively, in the triangulation corresponding to $\Gamma$.
In section~\ref{sec:triangles-graphs} we have shown
that $V=N_{lht}=$ the number of left-hand-turn paths in $\Gamma$,
$E=N_e=$ the number of edges in $\Gamma$,
and obviously, $F=N_v=$ the number of vertices in $\Gamma$.
Since $\Gamma$ is a 3-regular graph, we also have
$3F=2E$.
Hence, the Euler characteristic of \SCG{} 
can totally be recovered from $\Gamma$ by
$$\chi(\SCG)=V-E+F=N_{lht}-N_v/2.$$
So, the genus of \SCG{} is:
\begin{equation}
\label{eqn:genus}
  g=(2-\chi)/2=1+\frac{N_v-2N_{lht}}{4}\ .
\end{equation}

We could also calculate more geometrically
as follows:
The Euler characteristic, $\chi$, of \SOG{} is $2-2g-N_{lht}$.
We now use Gauss-Bonnet formula for
Riemannian surfaces without boundary :
$$\iint \kappa\,\dif\mbox{(area)} = 2\pi \chi,$$
where $\kappa$ is the curvature. In our case $\kappa \equiv -1$,
and the area of \SOG{} is $N_v\pi$,
since the area of the ideal triangle is $\pi$.
Substituing, we get
$-N_v\pi=2\pi(2-2g-N_{lht})$,
which leads again to formula~(\ref{eqn:genus}).

  %--------------------------------------------------
\section{Automorphisms}
\label{sec:automorphisms}
%---------------------------------------------------
Let $\Gamma_1$, $\Gamma_2$ be two graphs with orientation.
Let $\phi:\Gamma_1\to \Gamma_2$ be an isomorphism of graphs.

\begin{definition}
We say that $\phi$ is \emph{orientation-preserving} if
for any three edges, $\{e_1, e_2, e_3\}$, in $\Gamma_1$
meeting at a vertex $v$,
the cyclic ordered triplet $(e_1, e_2, e_3)$ fits
the orientation at $v$ iff
$\left(\phi(e_1), \phi(e_2), \phi(e_3)\right)$
fits the orientation at $\phi(v)$.
\end{definition}
A similar definition holds for \emph{orientation-reversing}.

We observe that the left-hand-turn paths are preserved under
an isomorphism of graphs with orientation.
In particular, the length of a left-hand-turn path is preserved.
These simple facts will help us later in finding automorphisms of
graphs.

In order to understand the automorphisms of a 3-regular
graph with orientation, let us first try to understand
the automorphisms of its universal cover, the 3-regular tree, $T_3$:
We choose some vertex of this tree, and call it the \emph{root}.
Denote by $\dot{T}_3$ the 3-regular rooted tree.
On a rooted tree we can define the notion of \emph{level}:
level $n$ is the set of all vertices which are at a distance $n$
from the root. A map between two rooted trees,
is a graph-map which maps the root
to the root. It has the property that it preserves the levels of the tree.

Now, Let $\dot{T}_3^{\mathrm{or}}$ be a 3-regular rooted tree, with
some choice of orientation on its vertices (actually, all the choices
are equivalent).
We have:
\begin{lemma}
  The group of automorphisms of $\dot{T}_3^{\mathrm{or}}$ is isomorphic
  to $\mathbb{Z}/3\mathbb{Z}$.
\end{lemma}
\begin{proof}
Since we are dealing with orientation-preserving automorphisms,
a generator is determined by the image of a first-level vertex of the tree.
\end{proof}

From the last lemma, we immediately get:
\begin{lemma}
  An automorphism of the
  3-regular tree with orientation, $T_3^{\mathrm{or}}$,
  is determined by the images of one vertex and an edge emanating from it.
\end{lemma}
Now, that we understand the automorphism group of the universal cover
we can proceed: 
\begin{lemma}
  Let $(v_0,e_0)$ and $(v_1,e_1)$ be two pairs of vertices and 
  edges emanating from them.
  Then, there exists at most one automorphism $\phi:\Gamma\to\Gamma$,
  such that $\phi(v_0,e_0)=(v_1,e_1)$.
\end{lemma}
\begin{proof}
  Let $\Gamma$ be a 3-regular graph with orientation,
  and $p:\tilde{\Gamma}\to\Gamma$ be its universal covering map
  (figure~\ref{fig:automorphisms}).
  \begin{figure}
    $$\includegraphics{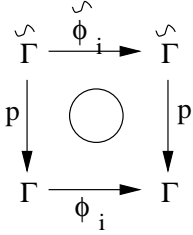}$$
    \caption{lifting of automorphisms}
    \label{fig:automorphisms}
  \end{figure}
  Choose liftings of $(v_i, e_i)$ $(i=0,1)$, $(\tilde{v}_i, \tilde{e}_i)$ 
  in $\tilde{\Gamma}$.
  If $\phi_1$ and $\phi_2$ both map $(v_0,e_0)$ to $(v_1,e_1)$ and 
  $\phi_1\neq\phi_2$ then they can be lifted to two
  different automorphisms 
  $\tilde{\phi}_1, \tilde{\phi}_2:\tilde{\Gamma}\to\tilde{\Gamma}$,
  which map $(\tilde{v}_0, \tilde{e}_0)$ to $(\tilde{v}_1, \tilde{e}_1)$.
  This contradicts the previous lemma.
\end{proof}
%-------------------------------------------------
The situation is very similar to what happens in $\D$:
We recall that an automorphism of the disk, $\D$ is determined
by the image of one point and a rotation.
This simple analogy shows:
\begin{proposition}
\label{prop:automorphisms}
Any automorphism of $\Gamma$ induces an automorphism
of \SOG.
\end{proposition}
%--------------------------------------------------
\subsection{Extension of Automorphisms}
Let $S$ be a complete hyperbolic Riemann surface of finite area,
and $i:S\hookrightarrow S_c$ be a conformal compactification 
(see section~\ref{sec:conformal}).
\begin{proposition}
  \label{prop:auto_extend}
  Any conformal automorphism of $S$ can be extended to a
  conformal automorphism of $S_c$.
\end{proposition}
\begin{proof}
Take a punctured-disk neighbourhood $U$ of a cusp $P$.
$i\circ\phi(U)$ is also a punctured disk.
Thus, by Riemann's theorem on removable singularities,
$P$ is a removeable singularity of $i\circ\phi$.
\end{proof}

  \section{Conformal Compactification}
\label{sec:conformal}
%-----------------------------------
Let $S$ be a Riemann surface, and
let $S_c$ be a compact Riemann surface.
If \mbox{$\phi:S\to S_c$} is a conformal embedding with dense image,
we say that $(S_c, \phi)$ is a \emph{conformal compactification} of $S$.\\
Examples:
\begin{enumerate}[\bfseries {Example} 1: \mdseries]
  \item
    The Riemann sphere is a conformal compactification of the complex plane.
  \item
      Take $S=\mathbb{H}^2$. By Riemann's Mapping Theorem
      $S$ is conformally equivalent to $\mathbb{C}-[1,\infty)$. Hence,
      the Riemann sphere is a conformal compactification of $\mathbb{H}^2$.
  \item \label{exmpl:finite_area}
     Suppose $S$ is a complete hyperbolic Riemann surface 
    of finite area.
    $S$ has a finite number of cusps
    (this follows from Gauss-Bonnet theorem or
    from Shimizu-Leutbecher theorem~\ref{thm:shimizu_leutbecher}).
    We may build a conformal compactification of $S$ as follows:
    We find pairwise disjoint neighbourhoods of the cusps,
    such that each neighbourhood is conformally equivalent to a
    punctured disk, and we fill in the missing point in each
    punctured disk. The conformal structure on the filled-in
    disk is unique, by the Uniformization Theorem 
    (see appendix~\ref{app:uniformization}).
    Thus, we get a compact Riemann surface 
    in which $S$ is conformally and densely embedded.
  \item  Any Riemann surface of infinite genus does
         not admit a conformal compactification.
%  \item  We suspect that the following is true:
%         A Riemann surface admits a conformal compactification
%         iff it admits a complete metric of finite total curvature.
%         (temp- maybe this is a theorem of A. Huber).
%        (temp- Is the following theorem true?:
%         \textbf{Theorem.}
%          S admits a complete Riemannian metric of finite area 
%          iff there do not exist bounded non-constant harmonic 
%          (analytic) functions on $S$).           
\end{enumerate}
%---------------------------------------------------------------
We will mainly be interested in example~\ref{exmpl:finite_area}.
\begin{proposition}
  In example~\ref{exmpl:finite_area} the conformal compactification
  is unique.
\end{proposition}
\begin{proof}
  Denote by $i:S\hookrightarrow S_{c,1}$ the conformal compactification 
  described above, and suppose that $\phi: S\hookrightarrow S_{c,2}$
  is a second conformal compactification.
  If the cusps of $S$ are $\{p_i\}_{i=1}^n$, then we get a conformal map
  $$\phi\circ i^{-1}:
  S_{c,1}-\{p_i\}_{i=1}^n\hookrightarrow S_{c,2}.$$
  The same arguments as in the proof of proposition~\ref{prop:auto_extend}
  show that $\phi\circ i^{-1}$
  can be extended to the cusps of $S$. So, we obtain that $S_{c,1}$ is
  conformally equivalent to $S_{c,2}$.
\end{proof}

  \section{Examples}
\label{sec:examples}
In order to give some examples, we first recall that
any \mbox{(anti-)auto}\-morphism of the graph $\Gamma$ induces an
\mbox{(anti-)auto}\-morphism of \SOG{}, and this, in turn
induces an \mbox{(anti-)auto}\-morphism of \SCG{}
(see section~\ref{sec:automorphisms}).
In other words, symmetries of the graph $\Gamma$ are reflected
in symmetries of \SCG{}.
% \SOG{} may have symmetries such as translation along geodesics
% which do not show in the graph.

Uniformize now \SCG{} by a metric of constant curvature.
The geodesic triangulation, $\mathcal{T}$, of $\SOG{}$
induced from $\Gamma$ maps
to some triangulation of $\SCG$ under the conformal embedding
$i:\SOG\hookrightarrow\SCG$. We may isotope $i(\mathcal{T})$ to
a geodesic triangulation of $\SCG$ keeping the vertices of
$\mathcal{T}$ fixed. Let us call this geodesic triangulation
$\mathcal{T}^{C}$.
The automorphisms of $\SCG{}$ which are extensions of automorphisms
of $\SOG$ stabilize $\mathcal{T}^{C}$.
%(temp- this is strange since they also stabilize the triangulation
%$\mathcal{T}$)
Thus, in order to determine the conformal structure of $\SCG{}$
we can follow the following steps:
\begin{enumerate}
  \item We determine the topological type of \SOG{}
        by counting the left-hand-turn paths in $\Gamma$ and using
        the genus formula~(\ref{eqn:genus}).
  \item We name the angles of the triangles in the geodesic triangulation,
        $\mathcal{T}^{C}$, of \SCG{}.
  \item For each left-hand-turn path in $\Gamma$ we write
        the corresponding cusp-equation: The angles around the filled-in 
        cusp amount to $2\pi$.
  \item \label{step:auto}
        We find the automorphisms and anti-automorphisms
        of the graph with orientation.
  \item We write the constraints on the angles implied by the symmetries
        found in step~\ref{step:auto}.
\end{enumerate}
Remark: In general these steps are not enough to determine
the conformal structure of $\SCG{}$, since we do not use
the full pasting rule in the construction of $\SOG{}$, which says 
that tick-marks go to tick-marks. 
We have yet to find a good method to take these conditions into account.

We begin by the simplest 3-regular graph, viz. two vertices and three edges.
It admits two possible non-isomorphic orientations:
\begin{enumerate}[\bfseries {Example} 1. \mdseries]
  \item
     See figure~\ref{fig:exmpl_v2g0}.
\begin{figure}
  $$\includegraphics{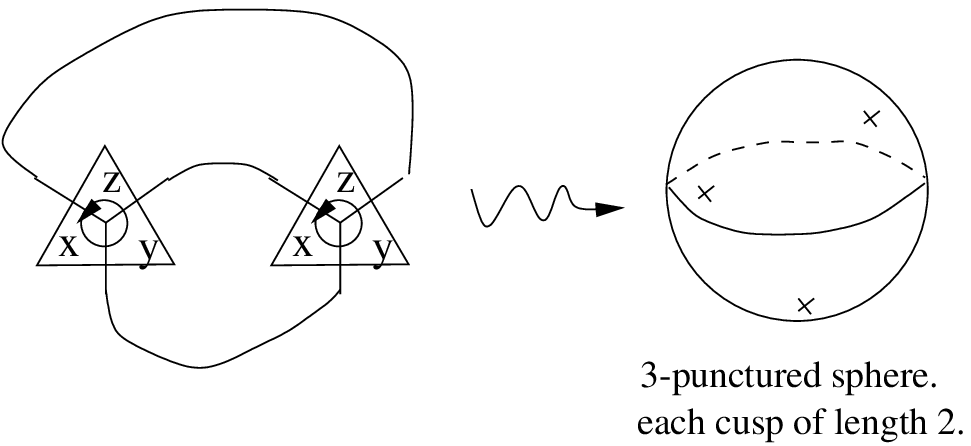}$$
  \caption{example no. \theenumi}
  \label{fig:exmpl_v2g0}
\end{figure}
We paste two triangles together to get a 3-punctured
sphere. Then we compactify to get the Riemann sphere
(the conformal structure on the sphere is unique).
     \label{exmpl:v2g0}
  \item
     See figure~\ref{fig:exmpl_v2g1}.
\begin{figure}
   $$\includegraphics{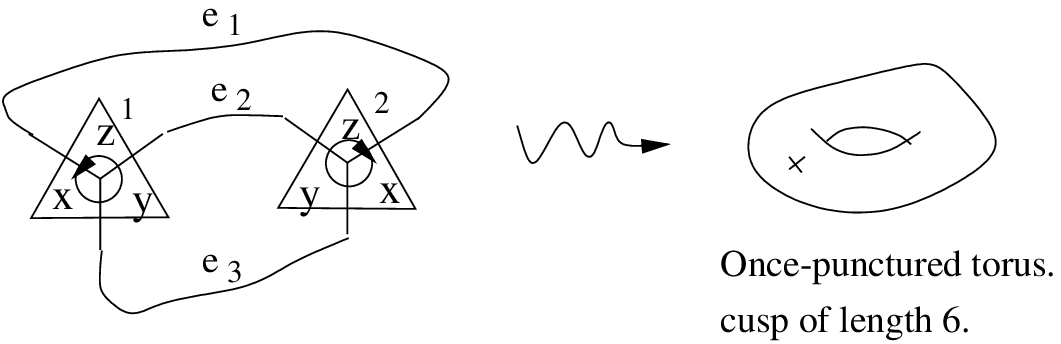}$$
   \caption{example no. \theenumi}
   \label{fig:exmpl_v2g1}
\end{figure}
We are still working with the graph of example~\ref{exmpl:v2g0}, 
but now we flip the right triangle over.
By the genus formula~(\ref{eqn:genus}),
we get a once-punctured torus.
We want to calculate the conformal class of this torus:
The angles around the filled-in cusp should amount to 2$\pi$:
\begin{equation*}
   x_1+z_2+y_1+x_2+z_1+y_2 = 2\pi,
\end{equation*}
We have the automorphism of order 6:
$$
\begin{array}{cc}
  v_1\mapsto v_2 &
  \begin{array}[t]{l}
  e_1 \mapsto e_2\\
  e_2 \mapsto e_3\\
  e_3 \mapsto e_1
  \end{array}
\end{array}
$$
which implies (together with its iterates) the equalities
$$ x_1=y_1=z_1=x_2=y_2=z_2. $$
Hence, we have that each angle equals
$\pi/3$, and we obtain the equilateral torus.

\end{enumerate}

In the next three examples we take $\Gamma$ to be the 1-skeleton of
the tetrahedron:
\begin{enumerate}[\bfseries {Example} 1. \mdseries]
  \setcounter{enumi}{2}
  \item
     See figure~\ref{fig:exmpl_ttrhdrn_f0}.
\begin{figure}
\includegraphics{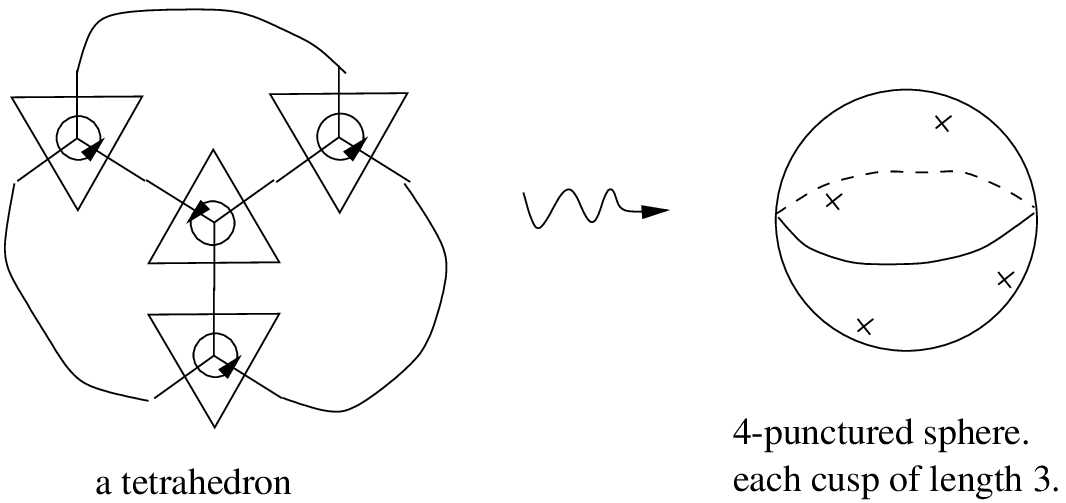}
\caption{example no. \theenumi}
\label{fig:exmpl_ttrhdrn_f0}
\end{figure}
We take the 1-skeleton of the tetrahedron,
with orientation induced from its surface.
The surface \SOG{} is a sphere
with 4 cusps at the vertices of a tetrahedron.
It is readily seen by the symmetries of the tetrahedron
that those vertices may be taken to be at
$0, 1, \omega, \omega^2$,
where $\omega$ is a primitive cube root of
unity. Obviously, \SCG{} is the Riemann sphere.
     \label{exmpl:ttrhdrn_f0}
  \item
     See figure~\ref{fig:exmpl_ttrhdrn_f1}.
\begin{figure}
  $$\includegraphics{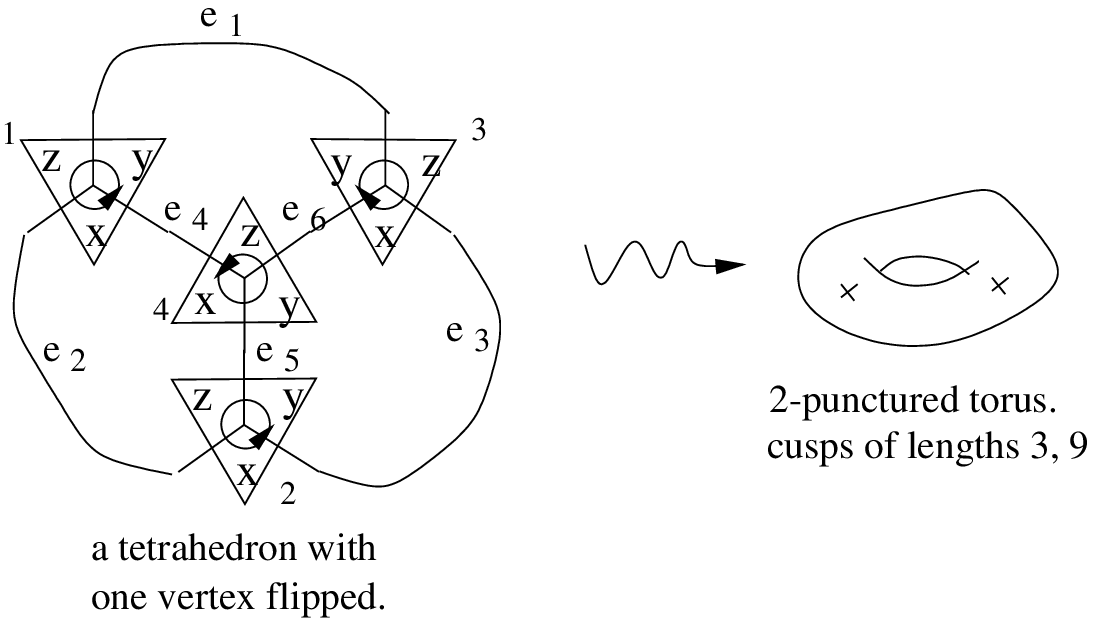}$$
  \caption{example no. \theenumi}
  \label{fig:exmpl_ttrhdrn_f1}
\end{figure}
We take the 1-skeleton of the tetrahedron
with the orientation induced from its surface,
and we flip orientation at one vertex.\\
\SOG\ is a 2-punctured torus.
The cusp equations are:
\begin{eqnarray*}
  x_1+z_2+x_4 &=& 2\pi \\
  y_1+z_4+x_3+x_2+z_1+y_3+y_4+y_2+z_3 &=& 2\pi.
\end{eqnarray*}
We have the automorphism of order 3 (rotation around $v_3$):
$$\begin{array}{cc}
    \begin{array}[t]{l}
    v_1\mapsto v_2\\
    v_2\mapsto v_4\\
    v_3\mapsto v_3\\
    v_4\mapsto v_1
    \end{array} &
    \begin{array}[t]{l}
    e_1\mapsto e_3\\
    e_2\mapsto e_5\\
    e_3\mapsto e_6\\
    e_4\mapsto e_2\\
    e_5\mapsto e_4\\
    e_6\mapsto e_1
    \end{array}
  \end{array}$$
which implies:
$$\begin{array}{l}
x_1=z_2=x_4,\\
y_1=x_2=y_4,\\ 
z_1=y_2=z_4,\\ 
x_3=y_3=z_3.
\end{array}$$
We have also an anti-automorphism of order 2:
$$\begin{array}{cc}
    \begin{array}[t]{l}
    v_1\mapsto v_1\\
    v_2\mapsto v_4\\
    v_3\mapsto v_3\\
    v_4\mapsto v_2
    \end{array} &
    \begin{array}[t]{l}
    e_1\mapsto e_1\\
    e_2\mapsto e_4\\
    e_3\mapsto e_6\\
    e_4\mapsto e_2\\
    e_5\mapsto e_5\\
    e_6\mapsto e_3
    \end{array}
  \end{array}$$
which implies:
$$\begin{array}{l}
y_1=z_1,\\
x_2=z_4,\\
y_2=y_4,\\
z_2=x_4,\\
y_3=z_3.
\end{array}$$
Combining all together, we have:
$$\begin{array}{l}
x_1=z_2=x_4=2\pi/3,\\
y_1=x_2=y_4=z_1=y_2=z_4=\alpha,\\
x_3=y_3=z_3=\beta,\\
6\alpha+3\beta=2\pi.
\end{array}$$
Since the triangulation on \SCG{} is geodesic
and our geometry is flat, we have the additional equations:
$$ x_i+y_i+z_i=\pi \;\; \mbox{for}\; 1\leq i\leq 4. $$
We obtain that triangle~3 is an equilateral triangle.
This implies that \SCG{} is an equilateral torus.

  \item
     See figure~\ref{fig:exmpl_ttrhdrn_f2}.
\begin{figure}
  $$\includegraphics{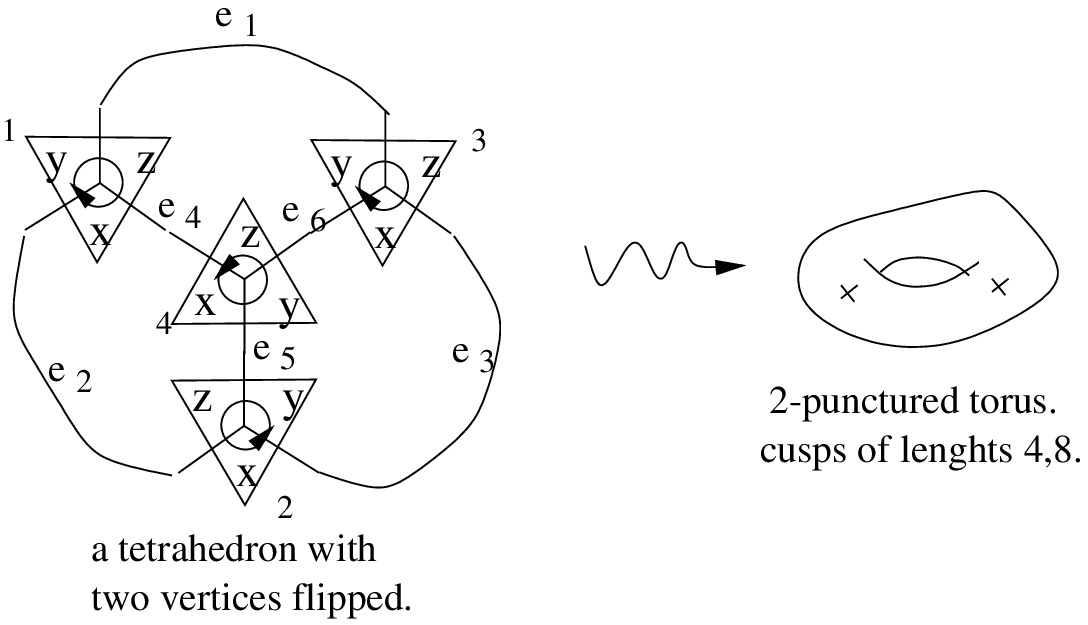}$$
  \caption{example no. \theenumi}
  \label{fig:exmpl_ttrhdrn_f2}
\end{figure}
We take the 1-skeleton of the tetrahedron
with the orientation induced from its surface,
and we flip orientation at two vertices.
$S^O(\Gamma)$ is a 2-punctured torus.
The cusp equations are:
\begin{eqnarray*}
  x_1+z_4+x_3+x_2 &=&2\pi,\\
  y_1+z_2+x_4+z_1+y_3+y_4+y_2+z_3 &=&2\pi.
\end{eqnarray*}
We have several anti-automorphisms of order 2:
$$
\begin{array}{cc}
  \begin{array}[t]{l}
  v_1\mapsto v_2\\
  v_2\mapsto v_1\\
  v_3\mapsto v_4\\
  v_4\mapsto v_3
  \end{array} &
  \begin{array}[t]{l}
  e_1\mapsto e_5\\
  e_2\mapsto e_2\\
  e_3\mapsto e_4\\
  e_4\mapsto e_3\\
  e_5\mapsto e_1\\
  e_6\mapsto e_6
  \end{array}
\end{array}
$$
which implies:
$$
\begin{array}{l}
  x_1=x_2,\\
  y_1=z_2,\\
  z_1=y_2,\\
  x_3=z_4,\\
  y_3=y_4,\\
  z_3=x_4.
\end{array}
$$
We have a second reflection:
$$
\begin{array}{cc}
  \begin{array}[t]{l}
  v_1\mapsto v_4\\
  v_2\mapsto v_3\\
  v_3\mapsto v_2\\
  v_4\mapsto v_1
  \end{array} &
  \begin{array}[t]{l}
  e_1\mapsto e_5\\
  e_2\mapsto e_6\\
  e_3\mapsto e_3\\
  e_4\mapsto e_4\\
  e_5\mapsto e_1\\
  e_6\mapsto e_2
  \end{array}
\end{array}
$$
which implies:
$$
\begin{array}{l}
  x_1=z_4,\\
  y_1=y_4,\\
  z_1=x_4,\\
  x_2=x_3,\\
  y_2=z_3,\\
  y_3=z_2.
\end{array}
$$
\ldots and a third reflection:
$$
\begin{array}{cc}
  \begin{array}[t]{l}
  v_1\mapsto v_3\\
  v_2\mapsto v_2\\
  v_3\mapsto v_1\\
  v_4\mapsto v_4
  \end{array} &
  \begin{array}[t]{l}
  e_1\mapsto e_1\\
  e_2\mapsto e_3\\
  e_3\mapsto e_2\\
  e_4\mapsto e_6\\
  e_5\mapsto e_5\\
  e_6\mapsto e_4
  \end{array}
\end{array}
$$
which implies:
$$
\begin{array}{l}
  y_1=z_3,\\
  z_1=y_3,\\
  x_1=x_3.
\end{array}
$$
combining all together, we have:
$$
\begin{array}{l}
  x_1=x_2=x_3=z_4 =\pi/2,\\
  y_1=z_1=y_2=z_2=y_3=z_3=x_4=y_4 =\pi/4.
\end{array}
$$
%(Question:Why do we get that the angles in the same left-hand-turn-path
% are all equal?)
Hence, we have a square torus.

\end{enumerate}

Next, we take $\Gamma$ to be the 1-skeleton of a cube:
\begin{enumerate}[\bfseries {Example} 1. \mdseries]
  \setcounter{enumi}{5}
  \item
     See figure~\ref{fig:exmpl_cube_f0}.
\begin{figure}
\includegraphics{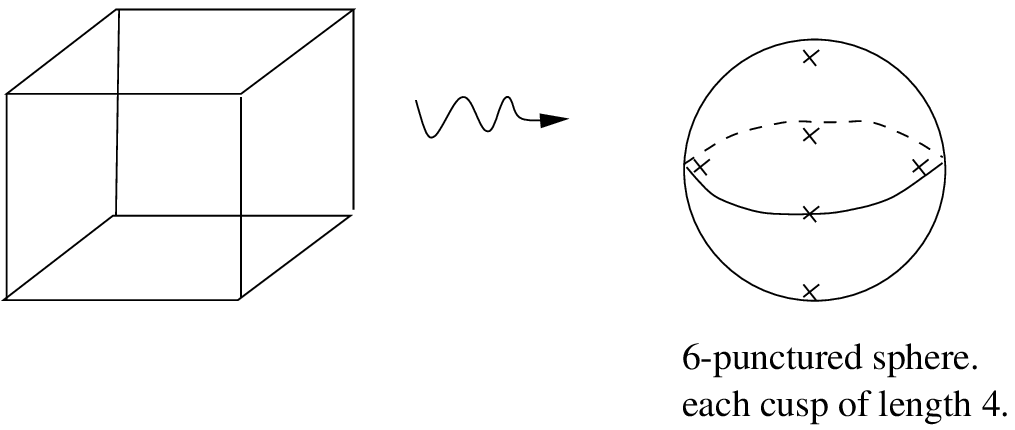}
\caption{example no. \theenumi}
\label{fig:exmpl_cube_f0}
\end{figure}
We take the 1-skeleton of the cube,
with orientation induced from its surface.
The surface \SOG{} is a sphere
with 6 cusps at the vertices of an octahedron.
It is readily seen by a stereographic projection that
the vertices can be taken to be $0,1,\mi,\infty,-1,-\mi$.
Then, we can map them by a M\"{o}bius transformation to
$1,\omega,\omega^2, -R, -R\omega,-R\omega^2$ respectively,
where $\omega$ is a primitive cube root of
unity and $R=2+\sqrt{3}$.
(If we apply another M\"{o}bius transformation: $z\mapsto -z/R$,
we get that we could take the vertices to be exactly at the vertices
and tick-marks of the standard ideal triangle).
\SCG{} is the Riemann sphere.

     \label{exmpl:cube_f0}
  \item
     We take the 1-skeleton of the cube
with orientation induced from its surface,
and we flip the orientation at one vertex, say $A$.
We calculate $\SCG$ according to ideas we learned from Curt McMullen.

Set $X=\SOG$. $X$ is a torus with four cusps.
One of the cusps is surrounded by twelve triangles and the other three are
surrounded by four each.
We have a $\mathbb{Z}/3\mathbb{Z}$ action on \SOG{},
coming from the order-3 rotation of the cube along its diagonal
through $A$.
The quotient, $Y$, is an orbifold with two cusps and
two ramification points of order 3:
One cusp is the image of the three cusps of length 4,
the other cusp is the image of the cusp of length 12,
and the two ramification points come from the centers of the
triangles pasted on the vertices of the rotation axis.\\
The genus of $Y$ may
be calculated by the Riemann--Hurwitz relation:
$$\chi(X)=3\chi(Y)-4,$$
where the relevant Euler characteristics are:
\begin{eqnarray*}
\chi(X)&=&2-2g_X-4,\\
\chi(Y)&=&2-2g_Y-2.
\end{eqnarray*}
After substituing we get that $g_Y=0$.

\noindent The compactification of $X$, $X^C$, is obtained as a 3-regular
covering of the sphere, branched over three points, each of them,
with ramification index 3, which may be taken to be $0$, $1$ and $\infty$.
Hence, a possible equation for $X^C$ in $\mathbb{C}\mathbb{P}^2$ is:
$$x^3=zy(y-z).$$

  \item
     (Due to Curt McMullen).
See figure~\ref{fig:exmpl_cube_fdiagonal}.
\begin{figure}
$$\includegraphics{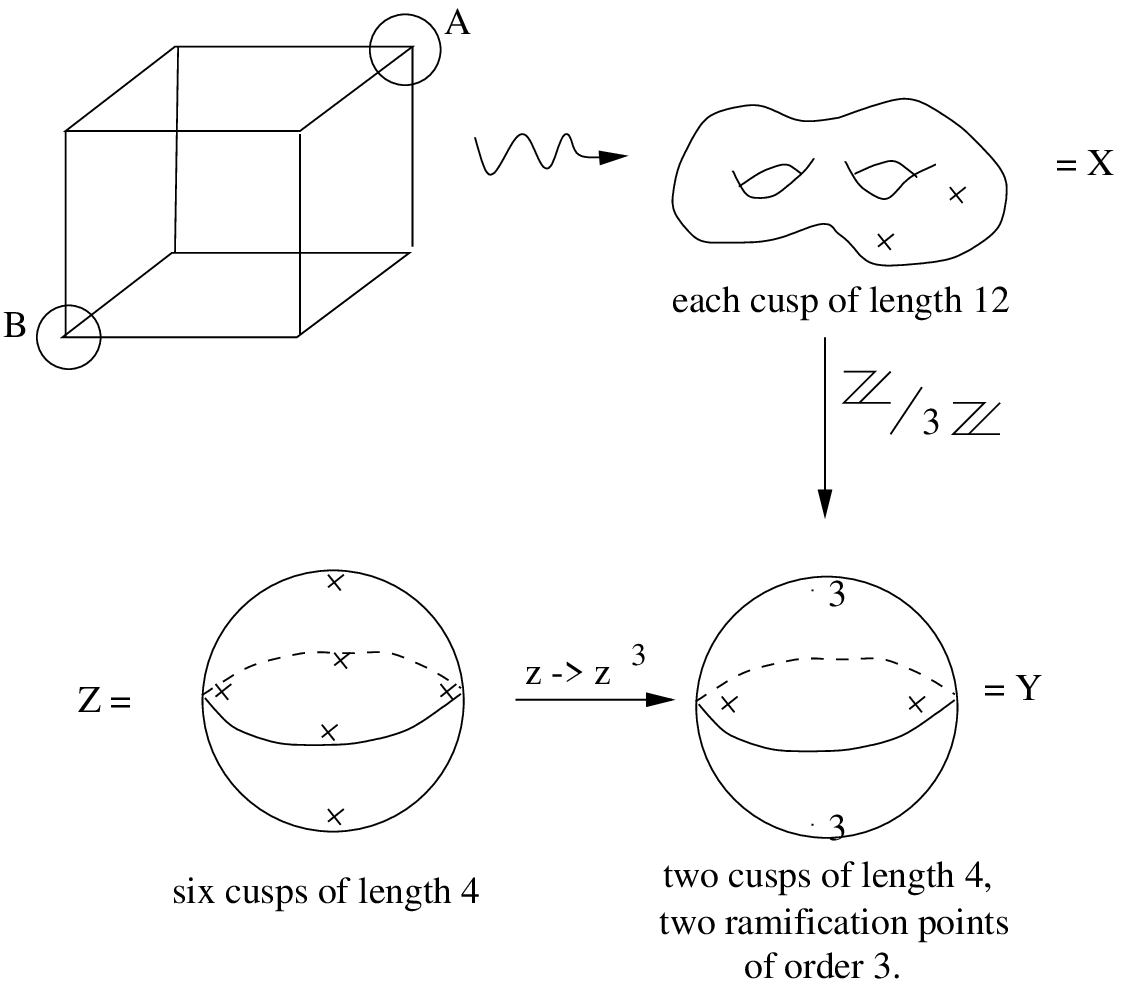}$$
\caption{example no. \theenumi}
\label{fig:exmpl_cube_fdiagonal}
\end{figure}
We take the 1-skeleton of the cube, with orientation
at the vertices induced from the surface of the cube.
Then we flip orientation at 2 opposite vertices, say $A$ \& $B$.
What we get is a surface, $X$, of genus 2, with 2 cusps,
and 12 ideal triangles around each cusp.

\noindent Consider the automorphism of $X$ of order 3 
obtained from rotating the cube along the diagonal $AB$.
When we mod out $X$ by the action of $\mathbb{Z}/3\mathbb{Z}$,
we get an orbifold, Y, with two cusps and two points
of ramification index 3:
The points of order 3 come from the centers of the triangles
at $A$ and $B$.
The genus of $Y$ may be calculated by the Riemann-Hurwitz relation:
\begin{eqnarray*}
  \chi(X) &=& 3\chi(Y)-4,\\
  \chi(X) &=& 2-2g_X-2 = -4,\\
  \chi(Y) &=& 2-2g_Y-2,
\end{eqnarray*}
from which we get $g_Y=0$.
The triangulation of $X$ descends to a triangulation of $Y$
composed of four cells around each cusp.
%(Question: Observe that if we mod out the graph, we no longer
%have a 3-regular graph - How can we see the quotient from the graph?)

\noindent Normalize $Y$ so that the ramification points
of order 3 are at $y=0$, and $y=\infty$,
and the cusps are at $y=1$, and $y=y_0$.
We now take the 3-regular cover, $Z$, of $Y$ with ramification
of order 3 at $y=0$ and $y=\infty$.
$Z$ is a surface (not orbifold) of genus 0, with 6 cusps
at the cube roots of $1$ and $y_0$.
One may verify this by again applying the Riemann-Hurwitz relation:
\begin{eqnarray*}
  \chi(Z) &=& 3\chi(Y)-4,\\
  \chi(Y) &=& 2-2g_Y-2 = 0,\\
  \chi(Z) &=& 2-2g_Z-6.\\
\end{eqnarray*}
The covering map $Z\rightarrow Y$ is given by $y=z^3$,
and the cellulation of $Y$ lifts to a triangulation of $Z$,
with four triangles surrounding each cusp.
Thus, it is readily seen that
$Z$ can be constructed from a 3-regular graph, $\Gamma$, which corresponds
to this triangulation.
$Z$ is of genus $0$ and has six cusps of length 4 each.
By formula~(\ref{eqn:genus}), $\Gamma$ has eight vertices.
\noindent The only graph with this combinatorics is the 1-skeleton
of the cube with orientation induced from the surface of the cube.
Therefore, $Z$ is a sphere with six cusps at the vertices of
an octahedron.
%Consider a sphere with cusps at the vertices of an octahedron,
%$\tilde{Z}$.
%$\mathbb{Z}/3\mathbb{Z}$ acts on $\tilde{Z}$,
%with quotient orbifold isomorphic to $Y$.
%Thus we have $Z=\tilde{Z}$.  

\noindent In example~\ref{exmpl:cube_f0}
it was shown that the vertices of an octahedron may be taken to be 
$\{1,\omega,\omega^2, -R, -\omega R,-\omega^2 R\}$, where $R=2+\sqrt{3}$
(or $R=2-\sqrt{3}$),
and $\omega$ is a primitive cube root of unity.
From here we conclude that $y_0=-R^3=-26-15\sqrt{3}$.
Hence, the conformal compactification of $X$, $X^C$,
is the surface obtained as a degree 3-regular
cover of $S^2$ with ramification points of order 3 
at $\{1, -26-15\sqrt{3}, 0, \infty\}$.
An equation for $X^C$ in $\mathbb{C}\mathbb{P}^2$ may be:
$$ x^3=\frac{zy(y-z)}{y+26z+15\sqrt{3}z}.$$

     \label{exmpl:cube_fdiagonal}
  \item
     (\cite{brooks:platonic, bfk}).
Set $\Gamma=PSL(2,\mathbb{Z})$, and
let $\rho$ be a finite index torsion free subgroup of
$\Gamma$.
We have:
\begin{theorem}
\label{thm:tor_free}
$\mathbb{H}^2/\rho$ can be constructed out of a 3-regular graph.
\end{theorem}
\begin{proof}
The classical fundamental domain for $\Gamma$ is given by 
(see figure~\ref{fig:psl2z_classic}):
$$F=\{z\in\mathbb{H}^2:|\Re(z)|<1/2, |z|>1\}.$$
\begin{figure}
$$\includegraphics{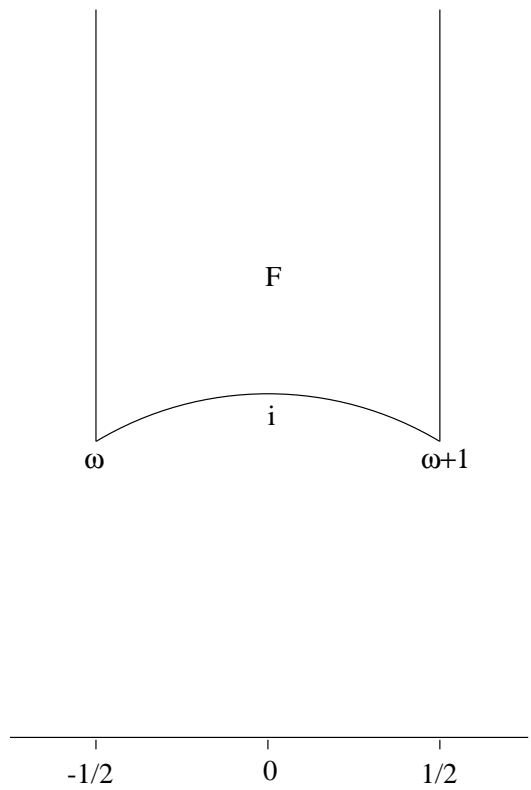}$$
\caption{Fundamental Domain for $PSL(2,\mathbb{Z})$}
\label{fig:psl2z_classic}
\end{figure}
Here, we take as a fundamental domain for $\Gamma$ 
(see figure~\ref{fig:psl2z}):
$$G=\{z\in\mathbb{H}^2:0<\Re(z)<1, |z| >1, |z-1| > 1\}.$$
\begin{figure}
  $$\scalebox{0.8}[0.8]{\includegraphics{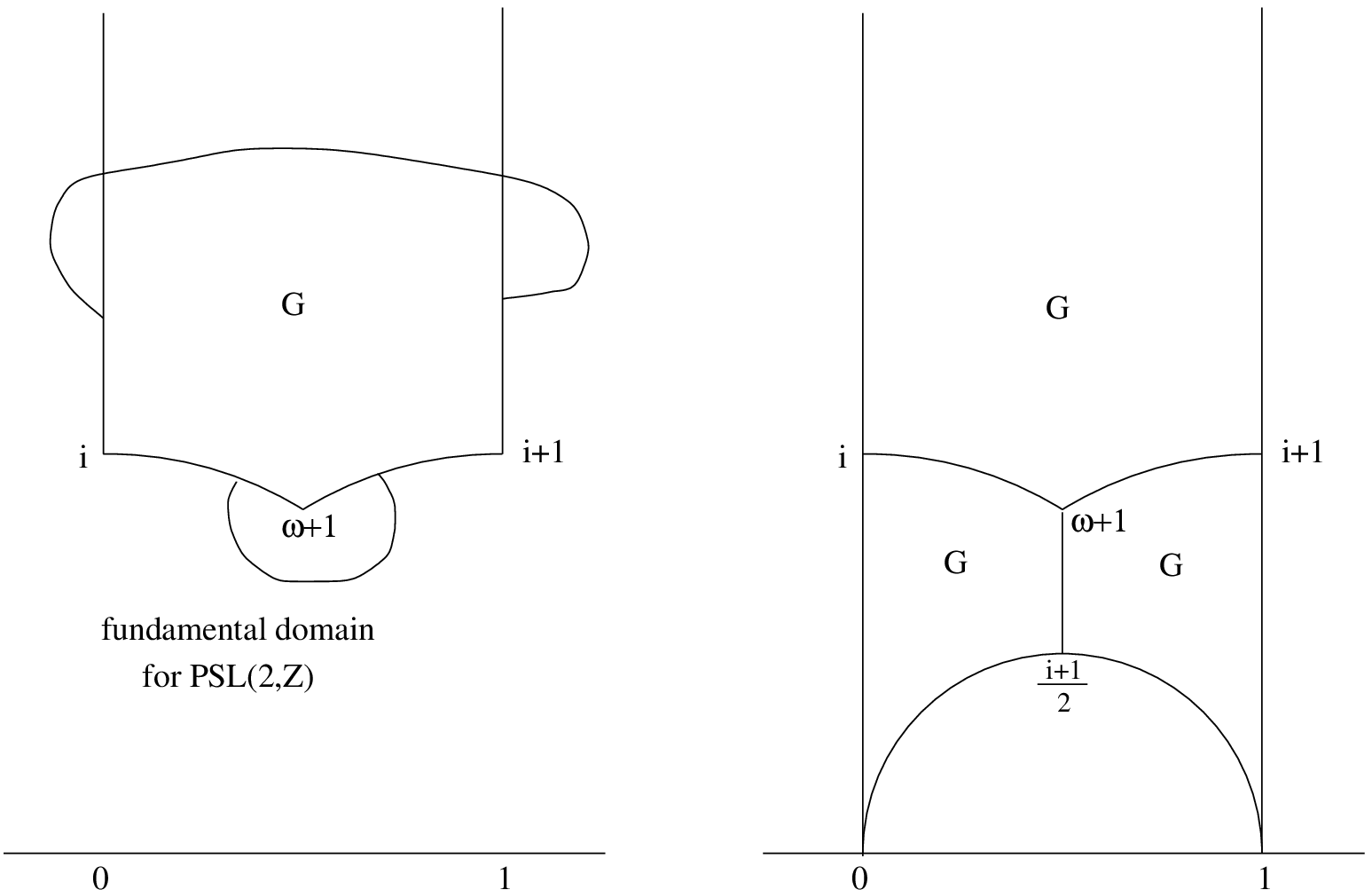}}$$
  \caption{3 copies of $G$ give the marked ideal triangle}
  \label{fig:psl2z}
\end{figure}
Observe that three copies of $G$ around $\omega + 1$ 
$(\omega=\me^{2\pi\mi/3})$
fit together to give the marked ideal triangle.
The equivalence between three such copies is given by an
elliptic element of order 3 in $\Gamma$.
 
A fundamental domain for $\rho$ is composed of copies
of $G$, and since $\rho$ is torsion free, the three
copies above are not equivalent under $\rho$,
and can all be included in a fundamental domain for $\rho$.
We have shown:
\begin{lemma}
A fundamental domain for $\rho$ can be chosen such that
it is composed of copies of the marked ideal triangle.
\end{lemma}
Let $J$ be such a fundamental domain for $\rho$, and mark
in each marked ideal triangle the three copies of $G$. 
Take $\Delta$ to be the three-regular graph which is 
composed of the finite edges of the copies of $G$, and its
vertices are the centers of the marked ideal triangles.
The orientaion on $\Delta$ is induced from the surface
$\mathbb{H}^2/\rho$. It is readily seen that
$\mathbb{H}^2/\rho \cong S^O(\Delta)$.
\end{proof}

A well known family of finite-index torsion free normal
subgroups of $\Gamma$ are the congruence subgroups 
(see~\cite{bfk}):
 $$\Gamma(k)=\!\!\left\{\left(
\begin{array}{cc}
  a & b \\
  c & d
\end{array}\right) \!\!\in \Gamma\ | 
\left(
\begin{array}{cc}
  a & b \\
  c & d
\end{array} \right)
\equiv
\pm\left(
\begin{array}{cc}
  1 & 0 \\
  0 & 1
\end{array}\right)
\!\!\!\!\!\pmod k \right\}.
$$
Observe that $\Gamma/\Gamma(k)\cong PSL(2,\mathbb{Z}/k\mathbb{Z})$.
As we have seen $S(k)=\mathbb{H}^2/\Gamma(k)$ all arise from graphs.
These graphs, called dual Platonic graphs in~\cite{bfk}, have
a simple algebraic description, and generalize the classical Platonic
solids.

Denote the conformal compactification of $S(k)$ by $P(k)$.
In \cite{brooks:platonic} these surfaces were called 
the \emph{Platonic Surfaces}.

     \label{exmpl:pltnc}
\end{enumerate}

  %-----------------------
\section{Belyi Surfaces}
%-----------------------
In this section we show that we know how to
construct Belyi surfaces out of graphs.
By Belyi's theorem it will follow
that the set of constructed surfaces is dense
in the moduli space.

Let $S$ be a compact Riemann surface.
It is well known that there exists a non-constant meromorphic
function on $S$, $\phi:S\rightarrow \mathbb{S}^2$.
\begin{definition}
If there exists a branched covering $\phi:S\rightarrow \mathbb{S}^2$,
such that $\phi$ is branched over at most three points,
then $S$ is called a \emph{Belyi surface}.
\end{definition}
\begin{theorem}[\cite{belyi}]
  $S$ is a Belyi surface if and only if as a curve in
  $\mathbb{C}\mathbb{P}^2$ its minimal polynomial lies
  over some number field.
\end{theorem}
\begin{corollary}
  For every integer $g\geq0$
  Belyi surfaces are dense in the
  moduli space of compact Riemann surfaces
  of genus $g$.
\end{corollary}
%--------------------------------------
We can characterize Belyi surfaces as follows:
\begin{lemma}
\label{lem:belyi}
  $S$ is a Belyi surface if and only if
  we can find finitely many points on $S$, $\{p_1,\ldots,p_k\}$,
  such that $S-\{p_1,\ldots,p_k\}$ is isomorphic to $\mathbb{H}^2/\rho$,
  where $\rho$ is a finite index torsion free subgroup of $PSL(2,\mathbb{Z})$.
\end{lemma}
\begin{proof}
  We begin by some basic observations:
  Let $\Gamma=PSL(2,\mathbb{Z})$,
  and $$\Gamma(k)=\left\{\left(
  \begin{array}{cc}
    a & b \\
    c & d
  \end{array}\right) \in \Gamma\ |
  \left(
  \begin{array}{cc}
    a & b \\
    c & d
  \end{array} \right)
\equiv
\pm\left(
\begin{array}{cc}
  1 & 0 \\
  0 & 1
\end{array}\right)
\pmod k \right\}.
$$
$\Gamma(2)$ has fundamental domain,F, which is composed
of 2 ideal triangles glued along a common edge (figure~\ref{fig:Gamma2}):
$$ F =\{ z\in \mathbb{H}^2 
: 0 < \Re(z) < 2, |z-1/2| > \frac{1}{2}, |z-3/2| > \frac{1}{2} \}. $$
\begin{figure}
  $$\includegraphics{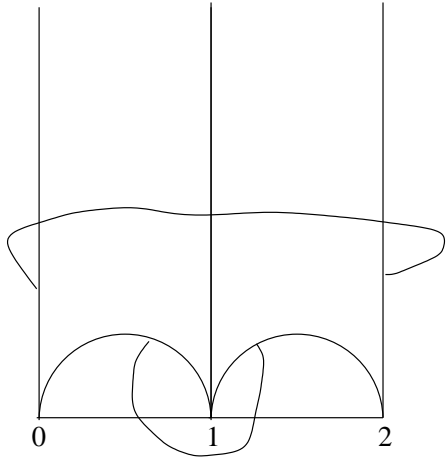}$$
  \caption{fundamental domain for $\Gamma(2)$}
  \label{fig:Gamma2}
\end{figure}
From this, we see that the 3-punctured sphere
(the punctures may be taken to be at $0, 1, \infty$)
may be realized as $\mathbb{H}^2/\Gamma(2)$.
Moreover, as we observed
in example~\ref{exmpl:pltnc}, each ideal triangle is
composed of three copies of the fundamental domain of $\Gamma$.
Therefore, the 3-punctured sphere is a degree-6 branched covering of
$\mathbb{H}^2/\Gamma$.

Let $S$ be a compact Riemann surface.
If $S$ is a Belyi surface, then there exists a finite number 
of points on $S$, $\{p_1,\ldots,p_k\}$ such that
$S - \{p_1,\ldots,p_k\}$ is a regular smooth finite degree covering of 
$\mathbb{H}^2/\Gamma(2)$.
Therefore, $S - \{p_1,\ldots,p_k\}$ may be realized as
$\mathbb{H}^2/\rho$, where $\rho$ is a finite-index 
subgroup of $\Gamma(2)$, and since $\Gamma(2)$ is a
finite-index torsion-free (see~\cite{bfk})
subgroup of $\Gamma$, we get that $\rho$ is a finite-index
torsion-free subgroup of $\Gamma$.

Conversely, if we can find finitely many points
$\{p_1,\ldots,p_k\}$ on $S$, such that $S - \{p_1,\ldots,p_k\}$
is isomorphic to $\mathbb{H}^2/\rho$, where $\rho$ is a finite index 
torsion free subgroup of $\Gamma$, then $S - \{p_1,\ldots,p_k\}$
is a finite-degree branched covering of $\mathbb{H}^2/\Gamma$, 
which is a sphere with one cusp and 2 ramification points.
Therefore, if we remove
the 2 ramification points and their pre-images, we get that
$S-\{p_1,\ldots,p_k,\ldots,p_n\}$ is a regular smooth 
finite-degree covering of the 3-punctured sphere. 
So, by definition, $S$ is a Belyi surface.
\end{proof}

Using the last lemma we can now prove:
\begin{theorem}
A Riemann surface can be constructed out
of a 3-regular graph if and only if it
is a Belyi surface.
\end{theorem}
\begin{proof}
Lemma~\ref{lem:belyi} and theorem~\ref{thm:tor_free}
show that any Belyi surface can be constructed 
out of some 3-regular graph with orientation.

Conversely, take a 3-regular graph with orientation, $\Delta$.
A fundamental domain, $F$, for $S^O(\Delta)$ is composed of copies
of the ideal triangle. By decomposing each ideal triangle into
three copies of the fundamental domain
for $PSL(2,\mathbb{Z})$, we see that $S^O(\Delta)$ may be realized as 
$\mathbb{H}^2/\rho$, where $\rho$ is a finite-index torsion 
free subgroup of $PSL(2,\mathbb{Z})$.
\end{proof}

%chapter 2
  \chapter{Geometric Relations between \SOG{} and \SCG{}.}
\label{ch:relations}
%---------------------
\section{Introduction}
%---------------------
Let $\Gamma$ be a 3-regular graph with orientation.
In chapter~\ref{ch:construction} we introduced a way
to construct a finite area Riemann surface, $\SOG$, out of $\Gamma$.
We denoted its conformal compactification by $\SCG$.

According to the Uniformization Theorem 
(see appendix~\ref{app:uniformization}),
\SCG{} admits a complete Riemannian metric
of constant curvature (which is compatible with the
given conformal structure).
In general, this metric might be very different from the metric
on \SOG{}. For instance, the 3-punctured
sphere admits a hyperbolic structure while its compactification,
the Riemann sphere, doesn't admit a hyperbolic metric at all.
  The main results of this chapter are theorems~\ref{thm:weak_control}
and~\ref{thm:comparison}. These results give conditions
under which we do have a hyperbolic metric on \SCG{} which 
is close in some sense to the metric on \SOG{}.
These conditions are related to the size of the cusps of \SOG{}.

As a last remark, we would like to suggest the reader
to bear in mind the \emph{Ahlfors--Schwarz Lemma} 
(see appendix~\ref{app:ahlfors_schwarz}) while reading this chapter.
%------------------
%-----------------------
\section{Large Cusps}
\label{sec:large_cusps}
%-----------------------
  The number of cusps on a complete hyperbolic Riemann surface of
finite area, $S$, is finite. This follows from Gauss--Bonnet theorem,
which reduces in that situation to:
$$ -\mbox{area}(S) = 2\pi\chi(S) = 2\pi(2-2g-N_c), $$
where $g$ is the genus of $S$, and $N_c$ is the number of cusps.
Alternatively, it follows from theorem~\ref{thm:shimizu_leutbecher}.

Each cusp of $S$ has \emph{a horocyclic neighbourhood}, i.e a neighbourhood
isometric to $\mathcal{B}_r=\left\{z\in\Dstar:|z|<r\right\}$ for
some $0<r<1$.
The area of $\mathcal{B}_r$, $\mathrm{area}(\mathcal{B}_r)$,
will be a measure for the size of the cusp:
\begin{definition}
  We say that \emph{the cusps of $S$ are bigger than $c$}
  if there exist pairwise disjoint horocyclic neighbourhoods
  of the cusps $\mathcal{U}_1,\ldots, \mathcal{U}_l$
  such that $\forall i\;\mathrm{area}(\mathcal{U}_i)>c$.
\end{definition}
%-------------------------

A beautiful classical theorem is the following:
\begin{theorem}[Shimizu-Leutbecher]
\label{thm:shimizu_leutbecher}
Let $S$ be a complete hyperbolic Riemann surface.
Then, the cusps of $S$ are $\geq$ 1.
\end{theorem}
The classical proof of Shimizu and Leutbecher uses an
iteration argument and can be found in~\cite{kra}.
We bring here another proof, based on the decomposition of $S$ 
into ideal triangles. A different proof in a similar spirit
  can be found in~\cite{buser}.
\begin{proof}
  We can find a locally finite triangulation on $S$,
  such that each vertex of the triangulation is a cusp of $S$.
  In each ideal triangle we can mark the 
  ``standard horocyclic segments'' 
  (see section~\ref{sec:ideal_triangles} 
  and figure~\ref{fig:triangle.poincare}). 
  Fix one cusp, $P$. Look at the ideal triangles
  which surround $P$. 
  We pick the triangle, $T$, with 
  the ``closest'' opposite horocycle to the cusp
  (in the universal cover, we can describe the situation
  easily - see figure~\ref{fig:shimizu_leutbecher}).
  \begin{figure}
    $$\scalebox{0.5}[0.5]{\includegraphics{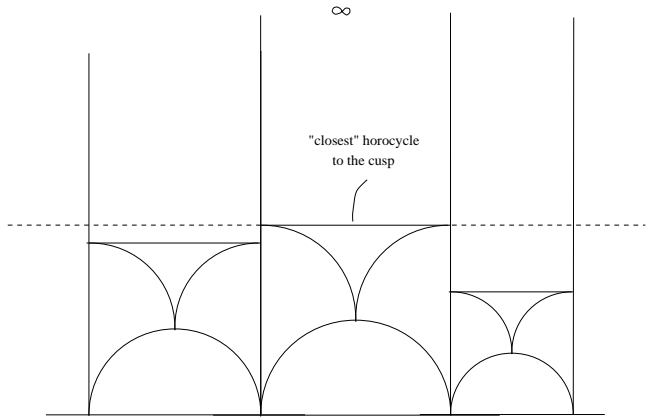}}$$
    \caption{proof of the Shimizu-Leutbecher Theorem}
    \label{fig:shimizu_leutbecher}
  \end{figure}
  We continue this horocycle segment to close a horocyclic
  neighbourhood around $P$.  
  We repeat this for all cusps of $S$. It is clear that
  the horocyclic neighbourhoods constructed in this way are pairwise
  disjoint, and since each horocycle segment is of length 1, the
  horocycles are of length at least~1.
\end{proof}

%-------------------------

%------------------
%-----------------------------
\section{Controlling $\SCG$}
%-----------------------------
Let $S$ be a complete hyperbolic Riemann surface
of finite area. Let $S_c$ be its conformal compactification
(see section~\ref{sec:conformal}).
The following theorem roughly says that
if the cusps of $S$ are large, we have some control
on the geometry of $S_c$.
\begin{theorem}
\label{thm:weak_control}
If the cusps of $S$ are bigger than $2\pi$,
then $S_c$ admits a conformal Riemannian metric
of negative curvature which coincides with
the metric on $S$ outside horocyclic cusp neighbourhoods,
$\{\mathcal{U}_i\}_{i=1}^l$, of area $>2\pi$.
\end{theorem}

In order to prove this theorem, we begin with some
preliminary calculations in the next two sections.
%---------------------
%--------------------------------------------------
\section{The Complete Hyperbolic Punctured Disk}
\label{sec:dstar}
%--------------------------------------------------
Let \Dstar{} denote punctured unit disk with
the Euclidean conformal structure on it.
The map from $\mathbb{H}^2$ onto \Dstar{} defined by
$z\mapsto\me^{2\pi\mi z}$ realizes \Dstar{} as
$\mathbb{H}^2/\{z\sim z+1\}$.
From this it is easy to give an explicit expression for
the complete hyperbolic metric on \Dstar{}:
\begin{lemma}
\label{lem:dstar_metric}
  The complete hyperbolic metric on \Dstar{}
  is given by
  $$ \dif s^2_{\Dstar} = \left(-\frac{1}{r\log r}\right)^2|\dif z|^2,$$
  where $r=|z|$.
\end{lemma}
\begin{proof}
  Set $w=\me^{2\pi\mi z}$. Then,
  \begin{eqnarray*}
  \dif s^2_{\Dstar} &=& 1/\Im(z)^2|\dif z|^2 \\
             &=& 1/\left(-\log |w|/2\pi\right)^2
                 \left|\frac{\dif w}{2\pi\mi w}\right|^2 \\
             &=& \left(-\frac{1}{|w|\log|w|}\right)^2|\dif w|^2.
  \end{eqnarray*}
\end{proof}

Denote by $\mathcal{B}_r$ the horocyclic neighbourhood
$\left\{z\in\Dstar:|z|<r\right\}$.
\begin{lemma}
  \label{lem:area}
  $\mathrm{area}(\mathcal{B}_r)=-2\pi/\log r$.
\end{lemma}
\begin{proof}
By lemma~\ref{lem:dstar_metric}
$$
  \begin{array}{lclclcl}
    \mathrm{area}(\mathcal{B}_r) &=&
      \displaystyle\int_0^{2\pi}\!\!\!\int_0^r -\frac{t}{(t\log t)^2}
                    \:\dif t\:\dif\theta
        &=& -2\pi\left.\displaystyle\frac{1}{\log t}\right|_0^r
        &=& \displaystyle{-\frac{2\pi}{\log r}.}
  \end{array}
$$
\end{proof}

We also have the following interpretation of the total geodesic
curvature of a closed curve:
\begin{lemma}
  \label{lem:geodesic_curvature}
  For a simple closed curve, $\gamma$,
  which bounds a domain $V$ which includes $0$, we have:
  $$\oint_\gamma \kappa_g\,\dif\mathrm{(length)} = \mathrm{area}(V).$$
\end{lemma}
\begin{proof}
  The Euler characteristic of $V$ is $0$.
  So, by the Gauss--Bonnet theorem we have:
  $$\oint_\gamma \kappa_g\,\dif\mathrm{(length)} +
    \iint_V \kappa\,\dif\mathrm{(area)} = 2\pi\chi(V)=0.$$
  But here $\kappa=-1$, and the lemma follows.
\end{proof}

\begin{figure}
  $$\includegraphics{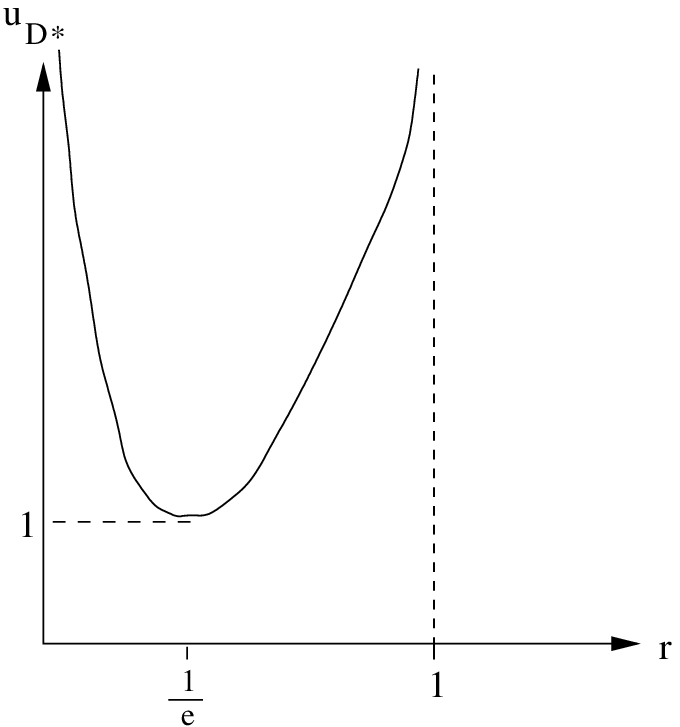}$$
  \caption{$\dif s_{\Dstar}^2=\me^{2u_{\Dstar}}(\dif x^2+\dif y^2)$}
  \label{fig:udstar}
\end{figure}

Set 
\begin{eqnarray}
  \lambda_{\Dstar} &=& -\frac{1}{r\log r},\\
  u_{\Dstar} &=& \log\lambda_{\Dstar}=\log\frac{1}{r\log\frac{1}{r}}.
  \label{eqn:udstar}
\end{eqnarray}
Then we have (see figure~\ref{fig:udstar}): 
\begin{lemma}
  \label{lem:udstar}
  $u_{\Dstar}$ satisfies:
  \begin{enumerate}[(i)]
    \item $\lim_{r\to 0} u_{\Dstar}(r) =\infty.$
    \item $\lim_{r\to 1} u_{\Dstar}(r) =\infty.$
    \item $u_{\Dstar}'(r)<0$ for $r<\frac{1}{\me}$.
    \item $u_{\Dstar}'(r)>0$ for $r>\frac{1}{\me}$.
    \item $\forall 0<r<1,\ u_{\Dstar}''(r)>0$.
  \end{enumerate}
\end{lemma}
\begin{proof}
Easy.
\end{proof}
%------------------------------------------------------------
\begin{comment}
\begin{lemma}
  \label{lem:per}
  $\mathrm{per}(\mathcal{B}_r)=
   \mathrm{area}(\mathcal{B}_r)=-2\pi/\log r$.
\end{lemma}
\begin{proof}
  By Lemma~\ref{lem:dstar_metric}
$$
  \begin{array}{lclclcl}
      \mathrm{per}(\mathcal{B}_r) &=&
      \displaystyle{\int_0^{2\pi} -\frac{1}{r\log r}r\dif\theta}
      &=& \displaystyle{-\frac{2\pi}{\log r}.}& & \\
    \mathrm{area}(\mathcal{B}_r) &=&
      \displaystyle\int_0^{2\pi}\!\!\!\int_0^r -\frac{t}{(t\log t)^2}
                    \:\dif t\:\dif\theta
      &=& -2\pi\left.\displaystyle\frac{1}{\log t}\right|_0^r
      &=& \displaystyle{-\frac{2\pi}{\log r}.}
  \end{array}
$$
In fact, one may prove that the area and the perimeter
are equal by Gauss--Bonnet theorem for $\mathcal{B}_r$:
Since $\partial \mathcal{B}_r$ is a horocycle,
its geodesic curvature is $1$.
Also, the Euler characteristic of $\mathcal{B}_r$ is $0$.
Therefore,
 $$ -\mbox{area}(\mathcal{B}_r) + \mbox{per}(\mathcal{B}_r) = 2\pi\cdot0,$$
leading to the desired result.
\end{proof}
\end{comment}
%---------------------------------------------------------
\section{The Curvature of Metrics on~\Dstar}
\label{sec:curvature}
%---------------------------------------------------------
Let $\dif s^2$ be a Riemannian metric on \Dstar{},
which is compatible with the conformal structure on $\Dstar$.
$\dif s^2$ is given by
\begin{equation}
  \label{eqn:radial_metric}
  \dif s^2 = \lambda(x,y)^2(\dif x^2+\dif y^2),
\end{equation}
where $\lambda$ is a positive function.
\begin{lemma}
\label{lem:curvature}
  The curvature of $\dif s^2$ is given by
  \begin{equation}
  \label{eqn:dstar_curvature}
    \kappa = -\frac{\Delta \log \lambda}{\lambda^2}.
  \end{equation}
\end{lemma}
\begin{proof}
  We calculate according to the following formula for
  the Gaussian curvature:
  \begin{equation}
  \label{eqn:curvform}
    \kappa=\frac{\left\langle\mathrm{K}(\dd{x},\dd{y})\dd{x},
                   \dd{y}\right\rangle}
           {\big\langle\dd{x},\dd{x}\big\rangle
            \big\langle\dd{y},\dd{y}\big\rangle},
  \end{equation}
  where
  \begin{equation}
    \mathrm{K}(X,Y) =\nabla_Y\nabla_X-\nabla_X\nabla_Y-\nabla_{[X,Y]}
  \end{equation}
  and $\nabla$ is the Levi-Civita connection. The Christoffel
  symbols for the Levi-Civita connection are given by
  \begin{equation}
  \label{eqn:chris}
    \Gamma_{ij}^{k}=\frac{1}{2}g^{km}(g_{im,j}+g_{mj,i}-g_{ji,m}),
  \end{equation}
  where $g_{ij}$ is the Riemannian metric tensor.
  Substituing the metric tensor $\dif s^2$ in this formula,
  we get\\
  \begin{minipage}{5cm}
  \begin{eqnarray*}
    \Gamma_{xx}^{x} &=& \frac{\partial\log \lambda}{\partial x} \\
    \Gamma_{xy}^{x} &=& \Gamma_{yx}^{x}= \frac{\partial\log \lambda}
                        {\partial y} \\
    \Gamma_{yy}^{x} &=& -\frac{\partial\log \lambda}{\partial x}
  \end{eqnarray*}
  \end{minipage}
  \begin{minipage}{5cm}
  \begin{eqnarray*}
    \Gamma_{xx}^{y} &=& -\frac{\partial\log\lambda}{\partial y} \\
    \Gamma_{xy}^{y}&=&\Gamma_{yx}^{y}= \frac{\partial\log\lambda}
                      {\partial x} \\
    \Gamma_{yy}^{y} &=& \frac{\partial\log\lambda}{\partial y} 
  \end{eqnarray*}
  \end{minipage}\hfill
  \begin{minipage}{1cm}   %produce equation number
    \begin{eqnarray}\end{eqnarray}
  \end{minipage}\\[2ex]
  We need also $\nabla_{\dd{x_i}}\dd{x_j}=\Gamma_{ij}^{k}\dd{x_k}$
  to calculate:
  \begin{eqnarray}
      \nabla_{\dd{x}}\dd{x} &=& \frac{\partial\log\lambda}{\partial x}\dd{x}-
                                \frac{\partial\log\lambda}{\partial y}\dd{y},
                                \nonumber\\
      \nabla_{\dd{x}}\dd{y} &=&
          \nabla_{\dd{y}}\dd{x}=\frac{\partial\log\lambda}{\partial y}\dd{x}+
                                \frac{\partial\log\lambda}{\partial x}\dd{y},
                                \label{eqn:nablas}\\
      \nabla_{\dd{y}}\dd{y} &=& \frac{\partial\log\lambda}{\partial x}\dd{x}+
                                \frac{\partial\log\lambda}{\partial y}\dd{y},
                                \nonumber
  \end{eqnarray}
  \ldots and\\
  \hspace*{\fill}
  \begin{minipage}{5cm}
  \begin{eqnarray*}
     \big\langle\dd{x},\dd{x}\big\rangle
     \big\langle\dd{y},\dd{y}\big\rangle &=& \lambda^4\\
     \big[\dd{x},\dd{y}\big] &=& 0.
  \end{eqnarray*}
  \end{minipage}\hfill
  \begin{minipage}{1cm}    %produce equation number
    \begin{eqnarray}\label{eqn:norms}\end{eqnarray}
  \end{minipage}
  \\[2ex]
  \noindent When we substitute equations~(\ref{eqn:nablas})
  --~(\ref{eqn:norms}) into the curvature formula~(\ref{eqn:curvform}),
  we get the desired result~(\ref{eqn:dstar_curvature}).
\end{proof}

In the special case of a radial metric on $\Dstar$, we have
\begin{corollary}
  \label{cor:radial_curvature}
  The curvature of
  $$ \dif s^2=\lambda^2(r)(\dif r^2+ r^2\dif\theta^2) $$
  is given by
    $$
    \kappa= -\frac{\displaystyle
         (\log\lambda)'' +
          \frac{1}{r}(\log\lambda)'}{\lambda^2}.
    $$
\end{corollary}
\begin{proof}
  The Laplacian in polar coordinates is given by
  $$
    \Delta=\frac{\partial^2}{\partial r^2}+
           \frac{1}{r}\frac{\partial}{\partial r} +
           \frac{1}{r^2}\frac{\partial^2}{\partial \theta^2}.
  $$
\end{proof}

%---------------------
%-----------------------
\section{A Key Lemma}
\label{sec:key_lemma}
%-----------------------
Recall that $\mathcal{B}_r =\{ z\in\Dstar :|z|<r\}.$
\begin{lemma}
  \label{lem:local}
  If $\mathrm{area}(\mathcal{B}_{r_0}) > 2\pi$, then there exists a Riemannian
  metric on $\dif s^2$ on \Dstar{}, such that
  \begin{enumerate}[i)]
    \item $\dif s^2$ is conformally equivalent to $\dif s^2_{\Dstar}$.
    \item $\dif s^2$ coincides with $\dif s^2_{\Dstar}$ outside a circle
          of Euclidean radius $r_0$.
          \label{con:geq_r0}
    \item $\dif s^2$ extends smoothly to a metric on $\D$.
    \item $\dif s^2$ has negative curvature.
  \end{enumerate}
\end{lemma}
\begin{proof}
By lemma~\ref{lem:area}, $\mathrm{area}(\mathcal{B}_{r_0})>2\pi$ implies
$r_0>\frac{1}{\me}$. By lemma~\ref{lem:udstar} $u_{\Dstar}'(r_0)>0$ and
$u_{\Dstar}''(r_0)>0$.  
Therefore, there exists a strictly convex monotonically increasing
$C^2$-function 
$u(r)$ on $[0,1)$ such that
$$ 
 \begin{array}{l}
  \forall r>r_0,\ u(r)=u_{\Dstar}(r),\\
  u'(0)=0.
\end{array}
$$
We can construct $u(r)$ as follows:
First, define a positive continuous function $w(r)$ such that,
\begin{eqnarray}
  \int_0^{r_0} w(r)\,\dif r & = & u_{\Dstar}'(r_0), \\ 
  \forall r\geq r_0,\ w(r)  & = & u_{\Dstar}''(r).
\end{eqnarray}
Next, define the function $v(r)$ by
\begin{equation}
  v(r)=\int_0^r w(s)\,\dif s. 
\end{equation}
Observe that $v(r)$ is a $C^1$-function which coincides with $u_{\Dstar}'(r)$
for $r\geq r_0$.
Finally, define $u(r)$ by
\begin{equation}
  u(r)=u_{\Dstar}(r_0)+\int_{r_0}^r v(s)\,\dif s.
\end{equation}
$u$ is a $C^2$-function which satisfies:
\begin{eqnarray}
  \forall r\geq r_0,\ u(r) &=& u_{\Dstar}(r),\\
  u'(0)&=&v(0)=0,\\
  \forall r>0,\ u'(r)&>&0\\
  u''=w>0 &\Rightarrow& \mbox{$u$ is strictly convex}.
\end{eqnarray}

Set $\lambda(r)=\me^{u(r)}$, and
$\dif s^2=\lambda^2(r)(\dif x^2+\dif y^2)$.
Then, we have:
\begin{enumerate}[a)]
  \item $\lambda'(0)=0 \Rightarrow$ $\dif s^2$ is smooth at $0$.
  \item $\dif s^2$ coincides with $\dif s^2_{\Dstar}$ for $r>r_0$.
  \item $\kappa(\dif s^2)<0$ by corollary~\ref{cor:radial_curvature}.
\end{enumerate}
\end{proof}
\section{Proof of Theorem~\ref{thm:weak_control}}
\label{sec:proof}
%---------------------------------------------------
The key lemma~\ref{lem:local} is the heart of the proof of
theorem~\ref{thm:weak_control}:
\begin{proof}
Let $\{\mathcal{U}_i\}_{i=1}^l$ be pairwise disjoint
horocyclic neighbourhoods of the cusps, $\{P_i\}_{i=1}^l$,
of areas $> 2\pi$.
According to the key lemma~\ref{lem:local}, for each $1\leq i\leq l$
we can alter the Riemannian metric in $\mathcal{U}_i$,
without affecting the metric outside of this neighbourhood.
Doing so, we are left with a Riemannian metric of negative curvature,
which extends smoothly across the cusps.
In other words, $S_c$ admits a Riemannian metric of negative curvature,
which coincides with the metric on $S$ outside the
horocyclic cusp-neighbourhoods, $\mathcal{U}_i's$.
\end{proof}

%------------------------------
%--------------------------------------------------------
\section{A Discussion of Theorem~\ref{thm:weak_control}}
\label{sec:discussion}
%--------------------------------------------------------
The best constant in theorem~\ref{thm:weak_control}
is $2\pi$.

One way to see this is
by applying the maximum principle for subharmonic functions:
Suppose we have a metric of negative curvature in $\D$ of the form
$$\dif s^2=e^{2u(x,y)}(\dif x^2+\dif y^2),$$
where $u$ is smooth and coincides with $u_{\Dstar}$ outside 
$\mathcal{B}_{r_0}$.
By formula~(\ref{eqn:dstar_curvature}) the negativity
of the curvature is equivalent to $\Delta u >0$.
Hence, $u$ is a subharmonic function.
By the maximum principle and Hopf's lemma
the maximum of $u$ in $\overline{\mathcal{B}_{r_0}}$
is attained on the boundary and $\frac{\partial u}{\partial r} >0$
on the boundary. But since $\frac{\partial u}{\partial r}$ on
the boundary = $u_{\Dstar}'(r_0)$, it follows by lemma~\ref{lem:udstar}
that $r_0>\frac{1}{\me}$. 
Finally, by lemma~\ref{lem:area} the last inequality is
equivalent to $\mathrm{area}(\mathcal{B}_{r_0})>2\pi$.

A different way to see that $2\pi$ is a best constant
is by applying the Gauss--Bonnet theorem 
to the region $\mathcal{B}_{r_0}$:
\begin{equation}
\label{eqn:gauss_bonnet}
  \iint_{\mathcal{B}_{r_0}} \kappa \,\dif\mbox{(area)} +
   \oint_{\partial\mathcal{B}_{r_0}} \kappa_g\,\dif\mbox{(length)} = 2\pi,
\end{equation}
where $\kappa$ is the curvature of the metric, and $\kappa_g$ is the
geodesic curvature of $\partial\mathcal{B}_{r_0}$.
Since $\kappa<0$ we obtain a necessary condition on
$\partial\mathcal{B}_{r_0}$:
\begin{equation}
\label{eqn:necessity}
  \oint_{\partial\mathcal{B}_{r_0}} \kappa_g\,\dif\mbox{(length)} > 2\pi.
\end{equation}
By lemma~\ref{lem:geodesic_curvature} this is equivalent to
$$\mathrm{area}(\mathcal{B}_{r_0})> 2\pi.$$
%--------------------------------
\section{An Extension Problem}
\label{sec:extension}
%--------------------------------
The necessary condition~(\ref{eqn:necessity}) raises a natural question:
Suppose we have a simple closed smooth curve $\gamma$ in $\Dstar$, 
which bounds a domain $V$ which contains $0$, such that 
\begin{equation}
\label{eqn:ext_necessity}
  \oint_{\gamma} \kappa_g\,\dif\mbox{(length)} > 2\pi.
\end{equation}
- Is it sufficient in order to extend the metric outside
of $V$ to a Riemannian metric with negative curvature in $V$?
A more general form of this question was raised already by M. Gromov 
in~\cite{gromov}, pp.~109--110.

Since the metric outside of $V$ is conformally equivalent to the
Euclidean metric, we may write:
$$\dif s_{\mathrm{out}}^2 = \lambda_{\mathrm{out}}^2(x,y)(\dif x^2+\dif y^2),$$
for some positive function $\lambda_{\mathrm{out}}$ defined outside of $V$.
We prove:
\begin{theorem}
\label{thm:no_extend}
If the maximum of $\lambda_{\mathrm{out}}$ on $\gamma$ is attained at
a point $z_{\mathrm{max}}$, where the normal derivative
$\frac{\partial \lambda}{\partial \vec{n}}(z_\mathrm{max})<0$, then
$\dif s^2$ cannot be extended to a conformal metric of negative curvature in
$V$.
\end{theorem}
\begin{proof}
Suppose $\dif s^2$ is a Riemannian metric of negative curvature 
on $\D$, which extends
$\dif s_{\mathrm{out}}^2$ and is given by 
$$\dif s^2=\lambda^2(x,y)(\dif x^2+\dif y^2),$$
where $\lambda$ is a smooth positive function defined on $\D$
and extends $\lambda_{\mathrm{out}}$.
By formula~\ref{eqn:dstar_curvature} 
$$\kappa(\dif s^2)<0\Rightarrow\Delta\lambda>0.$$
So, $\lambda$ is a subharmonic function.
By the maximum principle, the maximum of $\lambda$ in $\bar{V}$ is
attained on $\partial V=\gamma$.
Hence, it is attained at $z_\mathrm{max}$.
But, this is not possible, since the normal derivative there is negative.
\end{proof}

The next theorem shows that there exists curves as in 
theorem~\ref{thm:no_extend}:
\begin{theorem}
\label{thm:no_extend_example}
There exists a curve $\gamma$ in $\Dstar$, bounding a domain 
$V$ which includes $0$ such that
  \begin{enumerate}[i)]
     \item $\label{nec} \oint_{\gamma} \kappa_g\,\dif\mathrm{(length)} > 2\pi.$
     \item $\gamma$ has the property of theorem~\ref{thm:no_extend}
           with respect to the complete hyperbolic metric on $\Dstar$.
  \end{enumerate}
\end{theorem}
\begin{proof}
We can take $\gamma$ to be a horocycle with a slit
(see figure~\ref{fig:horocycle_slit}):
Smooth the curve whose image is given by
  $$  \{R_1e^{\mi t} : |t|<\theta\}
   \cup \{ R_2e^{\mi t}: \theta<|t|< \pi\}
   \cup\{re^{\mi t}: R_1\leq r\leq R_2, |t|=\theta\}. $$
  \begin{figure}
    $$\includegraphics{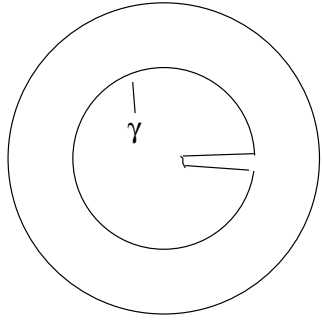}$$
    \caption{horocycle with a slit}
    \label{fig:horocycle_slit}
  \end{figure}
  We choose $R_2$ to satisfy $\frac{1}{\me}<R_2<1$. 
  By lemma~\ref{lem:udstar}, we can find $0<R_1<\frac{1}{\me}$ such that 
  $u_{\Dstar}(R_1)>u_{\Dstar}(R_2)$.
  Since by lemma~\ref{lem:area}
  $\mathrm{area}(\mathcal{B}_{R_2}) > 2\pi$,
  we can find $0<\theta<\pi$ such that the area bounded by $\gamma$
  in the complete hyperbolic metric on $\Dstar$ $>2\pi$.
  By lemma~\ref{lem:geodesic_curvature}, $\gamma$ 
  satisfies~\ref{nec}).

  By construction of $\gamma$ and by lemma~\ref{lem:udstar},
  we see that the maximum of $u_{\Dstar}=\log\lambda_{\Dstar}$ 
  on $\gamma$ is attained at a point where the normal derivative is
  negative. 
\end{proof}

If we work a little harder, we can show:
\begin{theorem}
\label{thm:no_extend_convex}
       $\gamma$ in theorem~\ref{thm:no_extend_example} can be taken 
       to be convex (i.e.\ $\kappa_g>0$).
\end{theorem}
\begin{proof}
Take $\gamma$ to be 
  be composed of a horocycle and a geodesic circular arc
  as in figure~\ref{fig:curve}.
  \begin{figure}
    $$\scalebox{.75}[.75]{\includegraphics{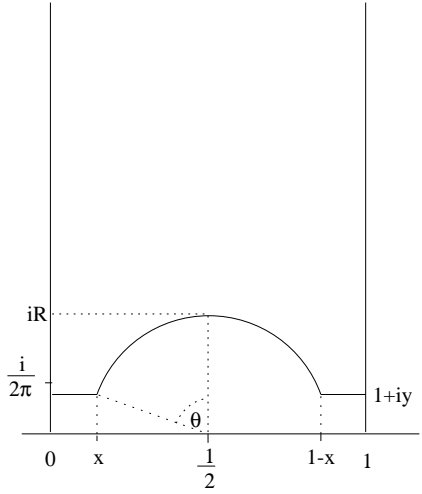}}$$
    \caption{The Extension Problem}
    \label{fig:curve}
  \end{figure}
  The total geodesic curvature (taking singular points into account)
  of $\gamma$ is given by the area above it 
  (lemma~\ref{lem:geodesic_curvature}):
  $$\oint_\gamma\kappa_g\,\dif\mbox{(length)} = \frac{2x}{y}+2\theta.$$
  We would like to have
  \begin{equation}
    \label{eqn:2pi}
    \frac{2x}{y}+2\theta>2\pi,
  \end{equation}
  and recalling that the map which realizes $\Dstar$ as 
  $\mathbb{H}^2/\{z\to z+1\}$ is given by $z\mapsto\me^{2\pi\mi z}$,
  we would also like to have
  \begin{equation}
    \label{eqn:max}
    u_{\Dstar}(\me^{-2\pi R(y,\theta)}) > u_{\Dstar}(\me^{-2\pi y}).
  \end{equation}
  We have
  \begin{eqnarray}
    x(y,\theta) &=& \frac{1-2y\tan\theta}{2},\\
    \label{eqn:R} R(y,\theta) &=& \frac{y}{\cos\theta}. 
  \end{eqnarray}
  Combining~(\ref{eqn:2pi})---(\ref{eqn:R}) together, and noting
  that $u_{\Dstar}$ is given by~(\ref{eqn:udstar}),
  we see that we should find $y>0, 0<\theta<\pi/2$ such that
  \begin{eqnarray}
    \label{eqn:2pi_explct}
    \frac{1}{y}-2\tan\theta+2\theta &>& 2\pi,\\
    \frac{2\pi y}{\cos\theta}-\log\frac{2\pi y}{\cos\theta}&>&
    2\pi y -\log 2\pi y.
  \end{eqnarray}
  We may first find a pair $(y_0,\theta_0)$ with equality 
  in~(\ref{eqn:2pi_explct}),
  and then decrease $\theta_0$. 
  So extracting $y$ from~(\ref{eqn:2pi_explct}), we obtain that we
  should find $0<\theta<\pi/2$ such that
  \begin{equation}
    \frac{\pi}{\pi+\tan\theta-\theta}>
    \frac{-\cos\theta\,\log\cos\theta}{1-\cos\theta}.
  \end{equation}
  The existence of such a $\theta$ can be shown by elementary calculus.
\end{proof}

Ultimately, we can show
\begin{theorem}
\label{thm:no_extend_all}
       There exists a convex curve $\gamma$ in $\Dstar$, 
       such that 
       \begin{enumerate}[i)]
         \item $\oint_{\gamma} \kappa_g\,\dif\mbox{(length)} > 2\pi$.
         \item no metric of negative curvature, conformal or
               nonconformal to the Euclidean metric, 
               can extend the outside metric 
               into the domain bounded by $\gamma$.
       \end{enumerate}
\end{theorem}
\begin{proof}
We construct $\gamma$ as follows (see figure~\ref{fig:no_extend_all}):
\begin{figure}
  $$\includegraphics{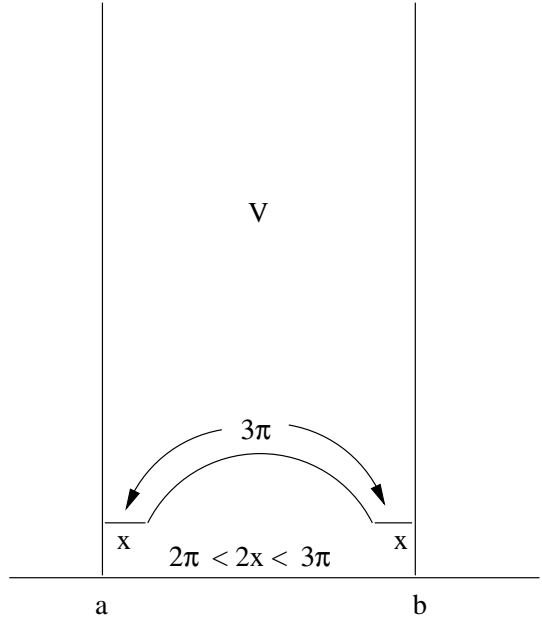}$$
  \caption{No extension}
  \label{fig:no_extend_all}
\end{figure}
Take any geodesic $\delta$ in $\mathbb{H}^2$, 
whose end points are in $\mathbb{R}$.
Take a symmetric segment $\delta'$ about the ``center'' of $\delta$ of length, 
say, $3\pi$. Now, extend the end points of $\delta'$, by horocyclic segments
of length $\pi$ each. Then, identify $\Dstar$ with 
$\mathbb{H}^2/\{z\mapsto z+C\}$, where $C$ is the Euclidean distance
between the endpoints of the horocyclic segments.

We get a convex curve $\gamma$ in $\Dstar$ whose total geodesic curvature
(taking singular points into account) $>2\pi$, by 
lemma~\ref{lem:geodesic_curvature}.
Suppose we have a Riemannian metric $\dif s^2$ of negative curvature
in the domain bounded by $\gamma$, $V$ ($\infty\in V$).
Since the distance between the end points of $\delta'$ $<3\pi$,
there exists a geodesic $\neq\delta'$ which connects these points.
But this is a contradiction to the fact that in a Riemannian surface of
negative curvature, there is only one geodesic in each homotopy class.
\end{proof}
\section{A Comparison Theorem}
\label{sec:comparison}
%------------------------------
If $S$ denotes a Riemann surface with a complete hyperbolic
metric of finite area $\dif s^2$,
and $S^C$ its conformal compactification,
we regard $S$ as embedded conformally in $S^C$.
The hyperbolic metric on $S^C$ will be denoted by $\dif {s^c}^2$.
%-------------------------------------------------------------
\begin{theorem}[\cite{brooks:platonic}]
  \label{thm:comparison}
  For every $\epsilon$, there exists $L(\epsilon)$ such that,
  if the cusps of $S$ $\geq L(\epsilon)$ then
  $$\frac{1}{1+\epsilon}\dif {s^c}^2\leq\dif s^2\leq 
  (1+\epsilon)\dif {s^c}^2.$$
\end{theorem}
%-------------------------------------------------------------
We will need the following lemma, which shows
that under a ``large cusps'' condition we can fix
the complete hyperbolic metric on $\Dstar$
to a smooth metric on $\D$ without losing control
on the curvature:
\begin{lemma}[\cite{brooks:platonic}]
  \label{lem:curvature_control}
  For any $\epsilon>0$ there exists $r_\epsilon$
  and a smooth metric $\dif s_\epsilon^2$ on $\D$,
  such that
  \begin{enumerate}[i)]
    \item
      $\forall r\geq r_\epsilon,\quad \dif s_\epsilon^2=
      \dif s_{\Dstar}^2$ ,\\
    \item
      $-(1+\epsilon) \leq \kappa(\dif s_\epsilon^2)\leq 
      \displaystyle\frac{-1}{1+\epsilon}.$
  \end{enumerate}
\end{lemma}
\begin{proof}
  We have that for a metric $\dif s^2=e^{2u}|\dif z|^2$,
  on the disk (or the punctured disk)
  the curvature is given by (see corollary~\ref{cor:radial_curvature}).
  \begin{equation}
    \kappa = -\frac{u''+\frac{1}{r}u'}{e^{2u}}.
  \end{equation}
  Set 
  \begin{equation}
    \label{eqn:g_of_u}
    g(r)\stackrel{\mathrm{def}}{=}\frac{\frac{1}{r}u'(r)}{e^{2u}}.
  \end{equation}
  In terms of $g(r)$, $u(r)$ is given by 
  \begin{equation}
    \label{eqn:u_from_g}
    e^{-2u(r)}=\int_r^1 2sg(s)\,\dif s.
  \end{equation}
  Hence the curvature is given in terms of $g(r)$ by
    $$\kappa=-2g(r)-rg'(r)-\frac{2r^2g(r)^2}{\int_r^1 2sg(s)\,\dif s}.$$
  From this formula we see that if we do not change $g(r)$ and
  $g'(r)$ too much then the curvature doesn't change much.
 
  We have already seen that the complete hyperbolic metrics on
  $\D$ and on $\Dstar$ are given, respectively, by
  \begin{eqnarray*}
     u_\D &=& \log\frac{2}{1-r^2},\\
     u_{\Dstar} &=& \log\frac{1}{r\log\frac{1}{r}}.
  \end{eqnarray*}
  One may calculate that (see \ref{eqn:g_of_u})
  \begin{eqnarray*}
    g_{\D}(r)  &=& \frac{1-r^2}{2},\\
    g_{\Dstar}(r) &=& \log{\frac{1}{r}}-(\log\frac{1}{r})^2.
  \end{eqnarray*}
  Hence,
  \begin{eqnarray}
    \lim_{r\to 1} g_\D(r) &=& \lim_{r\to 1} g_{\Dstar}(r) = 0,\\
    \lim_{r\to 1} g_\D'(r) &=& \lim_{r\to 1} g_{\Dstar}'(r) = -1
  \end{eqnarray}

  Thus, if $\delta>0$ is small,
  set $r_\epsilon$ such that
    $\forall r>r_\epsilon$
    \begin{eqnarray*}
       |g_\D(r)-g_{\Dstar}(r)|&<&\delta,\\
       |g_\D'(r)-g_{\Dstar}'(r)|&<&\delta,\\
    \end{eqnarray*}
  Then, define the function $\tilde{g}_\epsilon$ as follows:
  $$ \tilde{g}_\epsilon(r)=
  \left\{\begin{array}{lcr}
    g_\D(r)-g_\D(r_\epsilon)+g_{\Dstar}(r_\epsilon) &, & 
      0\leq r<r_\epsilon\\
    g_{\Dstar}(r)&,&
      r_\epsilon\leq r\leq 1.
   \end{array}\right.$$
  Finally, smooth $\tilde{g}_\epsilon$ at $r_\epsilon$
  to a smooth function, $g_\epsilon(r)$, 
  without altering much the derivative.

  For any $0\leq r\leq 1$, $g_\epsilon(r)$ and its
  derivative are close to $g_{\D}(r)$ or $g_{\Dstar}(r)$ and
  their derivatives. Therefore, the corresponding
  curvature $\kappa_\epsilon$ is close to $-1$.
  
  If we define $u_\epsilon$ by formula~(\ref{eqn:u_from_g}),
  then $u_\epsilon$ is smooth at $0$, it
  coincides with $u_{\Dstar}(r)$ for $r>r_\epsilon$ and
  its curvature is close to -1.
\end{proof}
%---------------------------------------------------------
Now we can prove the main theorem~(\ref{thm:comparison}).
\begin{proof}
We use the following version of the Ahlfors--Schwarz Lemma
due to Wolpert (see theorem~\ref{thm:wolpert_schwarz}):
%------------------------------------------------
\begin{theorem}[\cite{wolpert:schwarz}]
  Let $\dif s_1^2, \dif s_2^2$ be two (conformal) metrics on a compact
  Riemann surface. If $\kappa(\dif s_1^2)\leq\kappa(\dif s_2^2)<0$,
  then $\dif s_1^2 \leq \dif s_2^2$.
\end{theorem}
%------------------------------------------------
For any $\epsilon>0$
let $r_\epsilon$ be as in lemma~(\ref{lem:curvature_control}).
If the cusps of $S$ $>$ $\mathrm{area}(\mathcal{B}_{r_\epsilon})$, 
then $S^C$ admits a metric $\dif s_1^2$ whose curvature satisfies:
$$ -(1+\epsilon)<\kappa(\dif s_1^2) < -\frac{1}{1+\epsilon}, $$
and it conicides with the metric on $S$ away from
horocyclic neighbourhoods of the cusps.
We observe that 
$\kappa(\alpha\dif s^2)=\frac{1}{\alpha}\kappa(\dif s^2)$ 
(lemma~\ref{lem:curvature}). So, we have:
$$ \kappa(\frac{1}{1+\epsilon}\dif {s^c}^2)<
   \kappa(\dif s^2_1)<\kappa((1+\epsilon)\dif {s^c}^2)<0.$$
By Wolpert's theorem we get:
$$\frac{1}{1+\epsilon}\dif {s^c}^2\leq\dif s_1^2 \leq 
(1+\epsilon)\dif {s^c}^2.$$
\end{proof}
\appendix
  %---------------------------------
\chapter{Ahlfors--Schwarz Lemma}
\label{app:ahlfors_schwarz}
%---------------------------------
\section{The Lemma of Schwarz}
%------------------------------
Ahlfors, in his book \cite{ahlfors:complex}, proves
the lemma of Schwarz:
\begin{lemma}
Let $f:\mathbb{D}\to\mathbb{D}$ be an analytic function,
with $f(0)=0$.
Then 
\begin{equation}
  \label{eqn:ahl_int}
  |f(z)|\leq |z|
\end{equation}
and 
\begin{equation}
  \label{eqn:ahl_dif}
  |f'(0)|\leq1.
\end{equation}
If $|f(z)|=z$ for some $z\neq0$ or $f'(0)=1$ then $f$ is a rotation.
\end{lemma}
The proof is based on the maximum principle.
\begin{comment}
\begin{proof}
  Define
  $$ 
  g(z)=\left\{
  \begin{array}{ll}
    \frac{f(z)}{z}, & z\neq 0\\
    \scriptstyle f'(0), &z=0 
  \end{array}\right.
  $$
  $g$ is an analytic function on $\mathbb{D}$,
  and
  \begin{equation}
    \label{eqn:g_ineq}
    |g(z)|\leq \frac{1}{|z|} \;\mbox{for every}\; z\in\Dstar.
  \end{equation}
  Let $0<r<1$.
  $|g|\big|_{|z|\leq r}$ attains its maximum on $|z|=r$.
  Therefore, there exists $z_0$ with $|z_0|=r$,
  such that
  \begin{equation}
    \forall z, \; |z|\leq r \rightarrow |g(z)|\leq g(z_0) 
        \stackrel{(\ref{eqn:g_ineq})}{\leq} \frac{1}{r}.
  \end{equation}  
  Letting $r$ tend to $1$, we get $|g(z)|\leq1$ for all $z\in\mathbb{D}$.
  If $|g(z_0)|=1$ for some $z_0$, then $|g|$ attains its
  maximum in an interior point, so $g$ is constant function
  with absolute value $1$. Hence, $f$ is a rotation.
\end{proof}
\end{comment}

Ahlfors remarks that the condition $f(0)=0$ should
be considered only as a normalization:
If we set
$$ g_{\alpha}(z) = \frac{\alpha-z}{1-\overline{\alpha}z} $$
for $\alpha\in\D$, then 
$g_{\alpha}$ is an order-2 automorphism of the disk switching
$0$ and $\alpha$.
For any analytic map $f:\D\to\D$ and $z_1\in\D$ define
$$h=g_{f(z_1)}\circ f\circ g_{z_1}.$$
$h$ maps $\D$ into $\D$ and fixes $0$.
Applying the Lemma of Schwarz to $h$ with $z=g_{z_1}(z_2)$ we
obtain:
\begin{equation}
  \forall z_1,z_2\in\D\quad,
  \left|\frac{f(z_1)-f(z_2)}{1-\overline{f(z_1)}f(z_2)}\right|
   \leq \left|\frac{z_1-z_2}{1-\overline{z_1}z_2}\right|
\end{equation}
and in a differential form
\begin{equation}
  \forall z\in\D,\quad\frac{|f'(z)|}{(1-|f(z)|^2)}
  \leq\frac{1}{1-|z|^2}.
\end{equation}
Recalling that the complete hyperbolic Riemannian metric 
on $\mathbb{D}$ is given by 
$$\dif s^2=\frac{4}{(1-|z|^2)^2}|\dif z|^2, $$
or in its integrated form:
$$\tanh \mathrm{dist}(z_1,z_2) = 
  \left|\frac{z_1-z_2}{1-\overline{z_1}z_2}\right|,$$
we get a first geometric interpretation of the Lemma of Schwarz:
\begin{lemma}[Schwarz--Pick]
\label{lem:scwarz-pick}
If $f:\mathbb{D}\to\mathbb{D}$ is analytic,
then $f$ is length-decreasing. i.e., 
for any curve, $\gamma$,
$$\mathrm{length}(f(\gamma)) \leq \mathrm{length}(\gamma).$$
If $f$ preserves the distance between one pair of distinct points,
then $f$ is an isometry.
\end{lemma}

We have the following important corollary:
\begin{corollary}
  \label{cor:schwarz}
  Let $f:S_1\to S_2$ be an analytic mapping between complete
  hyperbolic Riemann surfaces. Then $f$ is length-decreasing.
\end{corollary}
\begin{proof}
  The universal covering of $S_i$ is $\D$ equipped with
  complete hyperbolic metric.
  We may lift $f$ to a an analytic map $\tilde{f}:\D\to\D$,
  between the universal coverings.
  $f$ is length-decreasing iff $\tilde{f}$ is length-decreasing.
  But $\tilde{f}$ is length-decreasing by the previous lemma. 
\end{proof}
%--------------------------------------------
\section{Curvature and the Lemma of Schwarz}
%--------------------------------------------
In \cite{ahlfors:schwarz} Ahlfors shows that the Lemma of Schwarz
is related with the notion of curvature.
\begin{lemma}[\cite{ahlfors:schwarz}]
Let $S$ be a Riemann surface endowed with a Riemannian metric $\dif s^2$
whose curvature $\kappa(\dif s^2)\leq -1$, and
let $f:\mathbb{D}\to S$ be analytic.
Then $f$ is length-decreasing.
\end{lemma}
\begin{proof}
We denote the complete hyperbolic metric on $\D$ by $\dif s_{\D}^2$.
Consider the pull-back singular metric on $\D$, $\dif s_f^2=f^*(\dif s^2)$.
$\dif s_f^2$ is conformally equivalent to $\dif s_{\D}^2$.
So, we may write:
\begin{eqnarray*}
   \dif s_{\D}^2 & = &\lambda_{\D}^2\left|\dif z\right|^2, \\
   \dif s_f^2    & = &\lambda_f^2\left|\dif z\right|^2,
\end{eqnarray*}
We need to show
\begin{equation}
  \label{eqn:ahlfors_schwarz}
  \lambda_f\leq\lambda_{\D}
\end{equation}
 (sometimes this is written as $\dif s_f\leq\dif s_{\D}$).
where the $\lambda$'s are nonnegative functions and $\lambda_{\D}>0$.
The last inequality is trivial at the branch points of $f$.

The curvatures (at nonsingular points) are given by (see~\ref{lem:curvature}):
  \begin{equation}
    \label{eqn:analytic_curvature}
    \kappa_i = -\frac{\Delta \log\lambda_i}{\lambda_i^2},
  \end{equation}
where $\Delta$ is the Laplacian: 
$\Delta=\partial^2/\partial x^2+\partial^2/\partial y^2$.
Hence, if we set $u_i=\log\lambda_i$ we have:
\begin{equation}
  \label{eqn:laplakappa}
  \Delta u_i=-\kappa_i\me^{2u_i}.
\end{equation}
Suppose first that $\frac{\lambda_f}{\lambda_{\D}}$ attains its maximum
in $\D$ (a maximum point is obviously nonsingular), 
then $u_f-u_{\D}$ attains its maximum in $\D$ and at maximum we have, 
\begin{eqnarray}
  \Delta (u_f-u_{\D})\leq0 &
  \stackrel{(\ref{eqn:laplakappa})}{\Rightarrow} &
   -\kappa_f\me^{2u_f}
   +\kappa_{\D}\me^{2u_{\D}}\leq0 \nonumber\\
   &\Rightarrow & 
    \me^{2(u_f-u_{\D})}
    \stackrel{(\kappa_{\D}<0)}\leq \kappa_{\D}/\kappa_f \leq 1 
          \label{eqn:ahlfors}\\
    &\Rightarrow& u_f-u_{\D}\leq 0 \nonumber\\
    &\Rightarrow& \frac{\lambda_f}{\lambda_{\D}}\leq1.
    \nonumber
  \end{eqnarray}
We got $\max\frac{\lambda_f}{\lambda_{\D}}\leq1$, as we wished.

In the general case, we pick $0<r<1$ and approximate $u_f$ by 
$$\lambda_{r,f}(z) = r\lambda_f(rz).$$
$\lambda_{r,f}$ is finite on $|z|=1$.
Hence $\frac{\lambda{r,f}}{\lambda{\D}}$
attains its maximum in $\D$. Also, $\kappa_{r,f}(z)=\kappa_f(rz)\leq -1$.
So, by~(\ref{eqn:ahlfors}) we get $\lambda_{r,f}\leq\lambda_{\D}$.
Letting $r$ tend to $1$ we obtain the desired result.
\end{proof}

One may ask what happens if we replace $\D$, with an arbitrary Riemann
surface.
For compact Riemann surfaces an answer was given by Wolpert:
\begin{theorem}[\cite{wolpert:schwarz}]
  \label{thm:wolpert_schwarz}
  Let $S_1$ and $S_2$ be two Riemann surfaces equipped
  with Riemannian metrics,
  where $S_1$ is compact and 
  $$\kappa_2\leq\kappa_1<0.$$
  Then, an analytic $f:S_1\to S_2$ is length decreasing. 
\end{theorem}
\begin{proof} 
  We pull-back the metric on $S_2$ to a metric on $S_1$.
  The two metrics on $S_1$ are conformal. Therefore their ratio is
  a well defined nonnegative function on $D$ which attains its maximum
  by compactness. From here, we proceed as in the proof
  of Ahlfors--Schwarz Lemma.
\end{proof}

Brooks proved the following version of Ahlfors--Schwarz Lemma:
\begin{theorem}[\cite{brooks:schwarz}]
  Let $S_1$ and $S_2$ be two Riemann surfaces equipped
  with Riemannian metrics and with $\D$ as an analytic univeral covering.   
  If $S_1$ is complete and 
  $$\sup\kappa_2\leq\inf\kappa_1<0, $$
  then an analytic $f:S_1\to S_2$ is length decreasing.
\end{theorem}
\begin{proof}
  We consider the univeral coverings, $\D_1$ and $\D_2$ with the pull-back
  metrics by the projection maps,
  and we lift $f$ to $\tilde{f}:\D_1\to\D_2$.
  Then we pull-back the metric on $\D_2$ to a metric on $\D_1$.
  The completeness of $S_1$ guarantees that the corresponding metric
  on $\D_1$ blows up on the boundary, and the $\sup-\inf$ inequality
  guarantees that the approximation argument in the proof
  of Ahlfors--Schwarz lemma will work.
\end{proof} 
%---------------------------------------------------
\begin{comment}
As an immediate corollary, we have:
\begin{lemma}
  Let $\dif s_1$ and $\dif s_2$ be two
  Riemannian metrics on a compact Riemann surface,
  such that there exist $C_1, C_2>0$ with
  $$C_1K(\dif s_2)\leq K(\dif s_1)\leq C_2 K(\dif s_2) < 0.$$
  Then $(1/C_1)\dif s_2^2\leq\dif s_1^2\leq(1/C_2)\dif s_2^2$.
\end{lemma}
\begin{proof}
  By the curvaure formula~(\ref{eqn:analytic_curvature})
  $K(\dif s/C) = C^2\cdot K(\dif s)$.
  Thus, we can apply Wolpert lemma. 
\end{proof}
\end{comment}
  \chapter{The Poincar\'{e} Polygon Theorem}
\label{app:PPT}
Let $\mathbb{H}^2$ be the hyperbolic plane.
Consider a finite-sided polygon
with all vertices at infinity, $P$, in $\mathbb{H}^2$.
$P$ has sides of two kinds:
The sides of the first kind are the sides which are not contained
in the boundary of $\mathbb{H}^2$.
The sides of the second kind are the sides, which are
contained in the boundary of the hyperbolic plane.
Suppose that the number of sides of the first kind is even,
and we divide them into pairwise disjoint pairs.
To each pair $\{s_i, s_j\}$ we associate an orientation-%
preserving transformation, $A_{ij}$, which maps $s_i$
onto $s_j$ in such a way that $A_{ij}(P)\cap P =\emptyset$.
We consider the subgroup, $G$, of $ISO^{+}(\mathbb{H}^2) \cong
PSL(2,\mathbb{R})$,
which is generated by the side-pairing transformations.

\begin{definition_ch}
A \emph{vertex-cycle transformation} is
an element of $G$ which stabilizes a vertex of the polygon $P$.
\end{definition_ch}

\noindent Remark: For any vertex $v$ of $P$
$\mathrm{Stab}_G(v)$ is a cyclic group.
\begin{theorem_ch}
If all the vertex-cycle transformations in $G$ are parabolic,
then $G$ is a discrete subgroup of $ISO^{+}(\mathbb{H}^2)$,
with $P$ as a fundamental domain.
\end{theorem_ch}

I liked De Rham's survey and proof of the theorem in
the compact and non-compact polygones cases in~\cite{ppt}.

  \chapter{The Uniformization Theorem}
\label{app:uniformization}

The Uniformization Theorem can be formulated as follows:
\begin{sloppy}
\begin{theorem_ch}
  Any simply connected Riemann surface is conformally
  equivalent to one and only one of the following surfaces:
\begin{enumerate}[(i)]
  \item The Riemann sphere.
  \item The complex plane.
  \item The unit disk.
\end{enumerate}
\end{theorem_ch}
\end{sloppy}
A proof can be found in~\cite{rs}.

\begin{corollary_ch}
  On a Riemann surface there exists a complete Riemannian metric of
  constant curvature 1, 0, or -1.
\end{corollary_ch}

\end{document}